\documentclass[11pt]{amsart}

\usepackage{amsmath,amsthm,amsfonts,amssymb, mathrsfs}
\usepackage[margin=1in]{geometry}
\usepackage[shortlabels]{enumitem}
\newtheorem{theorem}{Theorem}
\newtheorem{lemma}{Lemma}
\newtheorem{corollary}{Corollary}
\newtheorem{proposition}{Proposition}
\newtheorem{remark}{Remark}

\calclayout
\numberwithin{equation}{section}

\usepackage{amsmath,amsthm,amsfonts,amssymb, mathrsfs}
\numberwithin{equation}{section}
\usepackage{color}

\newcommand{\Var}{\text{Var}}
\numberwithin{equation}{section}
\newcommand{\be}{\begin{equation}}
	\newcommand{\ee}{\end{equation}}
\newcommand{\bes}{\begin{equation*}}
	\newcommand{\ees}{\end{equation*}}

\newcommand{\bP}{\mathbb{P}}

\newcommand{\R}{\mathbf{R}}
\newcommand{\RR}{\mathcal{R}}

\newcommand{\Z}{\mathbf{Z}}

\newcommand{\E}{\mathbb{E}}

\newcommand{\bean}{\begin{eqnarray*}}
	\newcommand{\eean}{\end{eqnarray*}}

\newcommand{\wu}{\widetilde{U}}

\usepackage{anysize}
 \usepackage{parskip}     

\title[Blowup for interacting SDEs]{Instantaneous blowup for interacting SDEs with superlinear drift}
\author{Mathew Joseph \and Shubham Ovhal}
\thanks{Research supported in part ANRF grant CRG/2023/002667 and a CPDA grant from the Indian Statistical Institute }

\address{Mathew Joseph\\ Statmath Unit\\ Indian Statistical Institute\\ 8th Mile Mysore Road\\ Bangalore 560059
} \email{m.joseph@isibang.ac.in}
\address{Shubham Ovhal\\ Statmath Unit\\ Indian Statistical Institute\\ 8th Mile Mysore Road\\ Bangalore 560059
} \email{rs\_math2104@isibang.ac.in}

\begin{document}
	
	\maketitle

\begin{abstract}%\\
We consider a system of interacting SDEs on the integer lattice with multiplicative noise and a drift satisfying the finite Osgood's condition. We show instantaneous everywhere blowup for initial profiles decaying slower than $\exp \left( -\sqrt{\big|\log |x|\big|}\right)$. We employ the splitting-up method to compare the interacting system to a one-dimensional SDE which blows up.


\end{abstract}

\section{Introduction}
We consider the following system of interacting stochastic differential equations (ISDEs):
\be \tag{ISDE} \label{eq:ISDE} dU_t(x) =\left(\mathscr{L} U_t\right)(x) dt + b\left(U_t(x)\right) dt + U_t(x) \,dB_t(x) ,\quad x \in \Z, \ee
with a bounded nonnegative initial profile $U_0$. Here $B_t(x)$ are independent one dimensional standard Brownian motions, %$\sigma: \R \to \R$ is a Lipschitz function satisfying 
%\be \label{eq:sigma}
%\LL  \le \left|\frac{\sigma(x)}{x}\right| \le \lip, 
%\ee
and $b:\R \to \R_+$ is a {\it locally} Lipschitz function that is {\it non-decreasing}, with $b(0)=0$ and satisfies {\it Osgood's condition}
\be \label{eq:osgood}
\int_\zeta^{\infty}\frac{dx}{b(x)} <\infty.
\ee
for some $\zeta>0$.
%To avoid pathological situations we will moreover assume that $b$ grows faster than linearly in the sense that for every $C>0$ there is a $K>0$ such that $b(x) \ge Cx$ for all $x\ge K$
Assume further that there is an $\eta >0$ such that for all large $x$ 
\be
\label{eq:b:cond2} \frac{b(x)}{x}\ge 1+\eta.
\ee
%Assume further that the function $f:(0, \infty) \to \R_+$ given by
%\be \label{eq:f}
%f(x):= \frac{b(x)}{x} 
%\ee
%is non-decreasing. The Osgood condition \eqref{eq:osgood} clearly implies that $\lim_{x\to \infty}f(x)=\infty$. 
%and satifies {\color{red} $\lim_{x \downarrow 0} f(x) > \frac65$ ?? See Proposition \ref{prop:compare} ??}. 

The operator $\mathscr L$ is the generator of a continuous time one dimensional symmetric random walk $\left(\mathcal {X}_t\right)_{t\ge 0}$ on $\Z$ with {\it finite range} $\mathcal{R}$, that is $\mathcal {X}_t = \sum_{j=1}^{N_t}\mathcal{Z}_i\,$, where $N_t$ is a Poisson process with rate one and $\mathcal{Z}_i$'s are i.i.d. $\Z$ valued, symmetric random variables satisfying $|\mathcal Z_1|\le \mathcal R$ almost surely, and $\Var (\mathcal{Z}_1) = \sigma^2$. These will be our standing assumptions in the rest of the paper and we skip repeating them in future.

{\it Mild} solutions of \eqref{eq:ISDE} are predictable processes which satisfy
\begin{align} \label{ISDE:int}
\begin{split}
U_t(x) & = \sum_{y \in \Z} G(t; y-x) U_0(y) + \sum_{y\in \Z} \int_0^t G(t-s; x, y ) \,b\left(U_s(y)\right)\, ds \\
&\qquad\qquad  + \sum_{y\in \Z} \int_0^t G(t-s; y-x ) \,U_s(y)\, dB_s(y) \hspace{2cm} \text{a.s.,}
\end{split}
\end{align}
where $G(t; y-x) := \bP\left(\mathcal{X}_t= y-x\right)$. As $\mathcal Z_1$ is symmetric we will sometimes write the first term above as $G_t*U_0(x)$. %A predictable process $t\mapsto U_t$ which solves the integral equation \eqref{ISDE:int} almost surely, and satisfies the moment bound
%\[
%    \E\left[ \int_0^t \sum_{y\in\Z} |U_s(y)|^2 ds \right] < \infty
%\]
%is called a {\it mild solution} to \eqref{eq:ISDE}. 
It is well known that if the function $b$ is Lipschitz continuous then the integral equation \eqref{eq:ISDE} has a unique continuous mild solution \cite{shig-shim,geor-jose-khos-shiu}, and if $b(0)=0$ then it is also nonnegative \cite{muel}.
%(we will sometimes use the notation $G(t; x):= G(t; 0,x)$), and 
%\[ \left(\widetilde G_t*U_0\right)(x)  = \sum_{y \in \Z} G(t; x, y) U_0(y).\]

%Let $(\Omega, \mathcal{F}, \mathcal{F}_t, \bP)$ be a filtered probability space such that $\{B_t(x)\}_{x\in\Z}$ are $\mathcal{F}_t$-measurable independent standard Brownian motions. If there exists a predictable process $U_t(x)$ which is an $\mathcal{F}_t$-measurable and almost surely continuous in $t$ at each site $x\in\Z$ such that $U_t(x)$ solves the integral equation \eqref{ISDE:int}, then the pair $\left(\{U_t(x)\}_{x\in\Z}, \{B_t(x)\}_{x\in\Z}\right)$ is called the a {\it mild} solution of \eqref{eq:ISDE}. If the function $b$ is Lipschitz continuous then the integral equation \eqref{eq:ISDE} has a unique continuous mild solution \cite{shig-shim}, and if $b(0)=0$ then it is also nonnegative \cite{muel}.}

For any fixed $J>0$,  consider the {\it truncated} system
\be \label{eq:UJ} dU_t^{(J)}(x) =\left(\mathscr{L} U_t^{(J)}\right)(x) dt + b\left(U_t^{(J)}(x) \wedge J\right) dt + U_t^{(J)}(x)\, dB_t(x). \ee
As $b$ is locally Lipschitz continuous, a mild solution $U^{(J)}_t(x)$ to equation \eqref{eq:UJ} exists; moreover since $b$ is non-decreasing, a comparison result \cite{geis-mant} shows that $U_t^{(J_1)}(x)\le U_t^{(J_2)}(x)$ when $J_1\le J_2$. As a consequence, for each $t>0$ and $x\in \Z$ the following limit exists almost surely
\be \label{soln}
U_t(x) := \lim_{J \to \infty} U_t^{(J)}(x).
\ee
We shall call the above the {\it minimal solution} to \eqref{eq:ISDE}.

Let us comment briefly on our notion of the minimal solution $U_t(x)$ in \eqref{soln}. First of all note that if $b$ is globally Lipschitz then it can be shown without much effort that the limit defined in \eqref{soln} is in fact the unique solution of \eqref{ISDE:int}.  Let $\mathcal{S}$ be the class of all non-negative globally Lipschitz continuous, non-decreasing functions. Let $f \in \mathcal{S}$ and $f \le b$, and suppose ${U}_{f}$ is the unique mild solution of equation \eqref{eq:ISDE} with $f$ in place of $b$. If ${U}^{(J)}_{f}$ are the corresponding truncated solution of ${U}_{f}$ then $ U^{(J)}_{f} \le U^{(J)}$ and consequently ${U}_{f} \le U$. %(See for instance \cite[Lemma 4.1]{geor-jose-khos-shiu}).
As $b(\cdot \wedge J) \in \mathcal{S}$ and $b(\cdot \wedge J) \le b $, it follows that
\[
    \sup_{f \in \mathcal{S}: f \le b } {U}_{f} = U.
\]
%by our assumptions on $b$. 
This clarifies the notion of `minimal solution' (see section 5 of \cite{foon-khos-nual}). %Note that if $b$ is globally Lipschitz then the limit defined in \eqref{soln} is in fact the unique solution of \eqref{ISDE:int}.% In summary if a `solution' to \eqref{eq:ISDE} is finite then it is at least larger than the limit defined in \eqref{soln} (see section 5 of \cite{foon-khos-nual}).

We will consider the following class of non-negative, symmetric and bounded initial profiles $U_0$ with a sufficiently slow decay to $0$. More precisely we assume in analogy with \cite{chen-khos-kim} the following: 
\begin{align}\label{u0:decay}
\begin{split}
&\text{ $U_0(x)$ is non-increasing on $[0,\infty)$, } \\
   % &U_0(x) = U_0(-x), \quad \quad U_0(0) \ge U_0(x), \\%\quad \text{ for all } x, \\
    &U_0(x) = U_0(-x),\qquad \lim_{|x|\to\infty} U_0(x) = 0, \qquad \lim_{|x|\to\infty} \frac{\left|\log U_0(x)\right|}{\sqrt{\log |x|}} = 0.
\end{split}
\end{align}
These conditions suffice to ensure instantaneous unboundedness in space for solution $U^{(0)}$ of \eqref{eq:UJ} with $J=0$ and is a crucial step in our method. Our main result is the following.

\begin{theorem}\label{thm:main} Consider \eqref{eq:ISDE} with the above assumptions on $b$ and $\mathscr L$. Let $U_t(x)$ be a minimal solution to \eqref{eq:ISDE} with initial profile $U_0$ satisfying \eqref{u0:decay}. Then we have instantaneous everywhere blowup of the system, that is for each $t>0$
\[ \bP\left(U_t(x) =\infty \text{ for all } x\in \Z\right) =1.\]
\end{theorem}

\begin{remark}
    By Mueller's comparison principle \cite{muel} and Theorem \ref{thm:main} we have the instantaneous explosion of \eqref{eq:ISDE} for any initial profile larger than those described in \eqref{u0:decay}, in particular for constant positive initial profiles.
\end{remark}

%\begin{remark} The above theorem can be extended to non-constant initial profiles. The proof of our theorem will in fact show the following. If the system without drift: $d\widetilde U_t(x) =\left(\mathscr{L} \widetilde U_t\right)(x)\, dt  + \widetilde U_t(x) \,dB_t(x)$ is unbounded for {\it all} $t>0$, that is
%\[ \bP \left(\sup_x \widetilde{U}_t(x) = \infty \right)=1 \text{ for all } t>0,\]
%then the solution to \eqref{eq:ISDE} instantaneously blows-up everywhere. Following \cite{chen-khos-kim} will then show that we have instantaneous blowup for \eqref{eq:ISDE} if the initial profile does not decay too rapidly. 
%\end{remark}

There has been recent interest in the blowup problems for stochastic heat equations (SHEs) with a drift term when the underlying space is $[0,1]$ or $\R$. This study started with \cite{bond-groi} where the authors show blowup (non-instantaneous) for the SHE $\partial_t u= \partial_x^2 u+b(u)+ \sigma\dot W$ on $[0,1]$ (with homogeneous boundary conditions) with {\it additive} white noise and a drift satisfying \eqref{eq:osgood}. Note that one cannot expect instantaneous blowup when the domain is bounded for bounded initial profile. Subsequently \cite{dala-khos-zhan} showed that \eqref{eq:osgood} is close to optimal for blowup for the SHE with drift in $[0,1]$. In particular they showed long time existence for drifts of the form $b(x)=O(|x||\log x|)$; they also considered multiplicative noise (noise of the form $\sigma(u) \dot W$). Later \cite{foon-nual} showed that \eqref{eq:osgood} is also a {\it necessary} condition for blowup for the SHE with {\it additive noise} on $[0,1]$; see also \cite{deng-liu-nane}.  More recently, \cite{foon-khos-nual} show the instantaneous everywhere blowup for solutions of the SHE $\partial_t u = \partial_x^2 u +b(u) + \sigma(u) \dot W$ on $\R$ with drift satisfying \eqref{eq:osgood} and $\sigma$ bounded away from $0$ and $\infty$.  As opposed to the handful of papers which show explosion for SHEs, there have been many more papers which prove global existence for SHE with superlinear drift and various types of noise, see for example \cite{chen-foon-huan-sali},  \cite{chen-huan}, \cite{sali}, \cite{sali-2}, \cite{sali-3}, \cite{shan-zhan}, \cite{shan-zhan-2} and the references therein.

The paper \cite{foon-khos-nual} on instantaneous blowup serves as a motivation for this paper. The important {\it unbounded} case of $\sigma(u)=u$ was not approachable by the methods in \cite{foon-khos-nual}. While our result is for interacting SDEs on $\Z$ (and not SHE on $\R$), it does indicate what to expect in the continuous case. The key idea in this paper is to apply the splitting-up method for SPDEs \cite{gyon-kryl} to deduce the blowup of the interacting SDEs from the blowup of the underlying SDE $dX_t = b(X_t)\, dt + X_t dB_t$. While the spliting-up method has been used extensively in the context of numerical approximations of SDEs and SPDEs, this is the first time (as far as we are aware) that it has been used for blowup problems of interacting SDEs. See also page 238 of \cite{muel} for an application to comparison results for SPDEs. 

There are three key ingredients of the proof of Theorem \ref{thm:main}. The first step is to show that for each arbitrarily small time $\delta>0$ the random field $U^{(0)}_{\delta}(\cdot)$ can take arbitrarily large values (see Proposition \ref{cor:sp:growth}); this follows from an application of the techniques from \cite{conu-jose-khos} and \cite{chen-khos-kim}. In particular for each large $M>0$ there is a point $p$ such that $U^{(0)}_{\delta}(p)\ge M$. The second key ingredient is Proposition \ref{prop:compare} where we show that the Alternating Process (the process obtained from the splitting method) started at $M\cdot \mathbf{1}\{x=p\}$, stays above the underlying SDE at site $p$ if $M$ is taken large enough. Hence the Alternating Process and consequently the ISDE becomes very large by time $2\delta$ at $p$. Finally we show that the $U^{(0)}$, started at this very large value at $p$, has to be large at an arbitrary site $x$ within the next time $\delta$ (see Proposition \ref{prop:det:noise}).

In Section \ref{sec:SDE} we discuss blowup for the underlying SDE. We define the Alternating Process for the truncated system \eqref{eq:UJ} in Section \ref{sec:AP} and show its process level convergence to \eqref{eq:UJ}. We also prove the key Proposition \ref{prop:compare} which compares the Alternating Process to the underlying SDE. Next we study the spatial growth of $U_t^{(0)}(\cdot)$ in Section \ref{sec:ISDE:no_drift} and compare the deterministic term of $U_t^{(0)}$ to the stochastic term for small time $t$. Finally we prove Theorem \ref{thm:main} in Section \ref{sec:thm}.
%\begin{itemize}
%\item {\color{red} Boinder Grossman paper}
%\item { \color{red} Notation: For a sequence $a_n, b_n$ we write  $a_n \lesssim b_n$ means $a_n \le C b_n$ for all $n$, for some universal constant $C$. $a_n\lesssim_{(M,J)} b_n$ means $a_n \le Cb_n$ where $C$ might depend on $M, J$. }
%\end{itemize}

\textbf{Notation:} For a sequence $a_n, b_n$ we write  $a_n \lesssim b_n$ to mean $a_n \le C b_n$ for all $n$, for some universal constant $C$. The notation $a_n\lesssim_{(M,J,\cdots)} b_n$ means $a_n \le Cb_n$ where $C$ might depend on $M, J,\cdots$. For any function $f$ on $\Z$ we denote it's supremum norm on $\Z$ by $\|f\|_\infty$. The maximum and minimum of two real numbers $a$ and $b$ will be denoted by $a\vee b$ and $a\wedge b$, respectively.

\section{The underlying SDE} \label{sec:SDE}

Let $B_t$ be a standard Brownian motion. For a locally Lipschitz continuous drift function $b$ consider the SDE \be \label{pro:X} dX_t= b(X_t) \,dt + X_t\, dB_t.\ee
%\begin{enumerate}
%\item {\color{red} Perhaps increasing convex ordering can be used for comparing general $\sigma(x)$ with $\sigma(x) = Cx$ ?} 
%\item {\color{red} Mention what is meant by explosion}
%\end{enumerate}
The main result of this section (Proposition \ref{prop:explosion}) is a variation on Osgood's Theorem on ODEs for this SDE. %, the main difference being a behavior dependent on initial value. 
In the following proposition we prove that the process $X_t$ takes on arbitrarily large values in a prescribed time with very high probability provided that the process $X_t$ starts with an initial value sufficiently large. 

The SDE \eqref{pro:X} which we consider here does not appear to be easily amenable to a straightforward analysis using the well known {\it Feller's explosion test} \cite[Theorem 5.29]{kara-shre}, therefore we proceed to study the SDE directly using first principles.

We briefly outline how a local in time solution of \eqref{pro:X} can be constructed when $b$ is non-negative, locally Lipschitz and increasing. By replacing the drift $b(x)$ in \eqref{pro:X} by the Lipschitz continuous function $b(x \wedge J)$ for any $J\ge0$ we obtain a corresponding SDE which admits a unique solution process $X_t^{(J)}$ with continuous paths. Define the stopping times $\theta_J := \inf\{t>0 : X_t^{(J)} \ge J\}$. Since $b$ is increasing we have $\theta_J\le \theta_{J+1}$ for all $J$. The explosion time is then defined to be $\theta_\infty := \lim_{J\to\infty} \theta_J$. The consistency relation
\[
X^{(J+1)}_t\mathbf{1}_{\{t \le \theta_{J}\}} = X^{(J)}_t\mathbf{1}_{\{t \le \theta_{J}\}}
\]
allows us to the define $X_t$ up to the explosion time $\theta_\infty$ as $X_t := \lim_{J\to\infty} X_t^{(J)}$. We denote by $\bP_K$ the conditional probability measure $\bP(X_t\in \cdot | X_0 = K)$.

For a fixed time $t>0$ we say that the process $X$ {\it explodes} to $\infty$ by time $t$ if $\theta_\infty<t$, or equivalently if there exists $s\le t$ such that $\lim_{r\uparrow s}X_r=\infty$.

\begin{proposition} \label{prop:explosion} Assume that the drift $b$ is locally Lipschitz, non-decreasing and the Osgood condition \eqref{eq:osgood} is satisfied. Let the process $X_t$ be the solution to the SDE \eqref{pro:X}. Then for every $\epsilon, \delta>0$ there exists $K=K(\epsilon, \delta)$ such that 
\[ \bP_{2^K}\left(X \text{ explodes to } \infty \text{ by time } \delta\right) \ge 1-\epsilon.\]
\end{proposition}

Before we proceed to the proof of the proposition we state two lemmas which will be required. Let 
\be \label{eq:f}
f(x):= \frac{b(x)}{4x} 
\ee

\begin{lemma} \label{lem:osgood} Under the hypothesis that $b$ is non-decreasing, the Osgood condition \eqref{eq:osgood} on $b$ is equivalent to 
\[ \sum_{k\ge 0} \frac{1}{f(2^k)} <\infty.\]
\end{lemma}
\begin{proof}
%Non-negativity of $b(x)$ allows us to write
   % \begin{align*}
      %  \int_1^\infty \frac{dx}{b(x)} = \sum_{k\ge0} \int_{2^k}^{2^{k+1}} \frac{dx}{b(x)}.
    %\end{align*}
    We partition the interval $[1, \infty)$ into the intervals $[2^k, 2^{k+1})$, and use the assumption that $b$ is non-decreasing to obtain 
    \begin{align*}
        \sum_{k\ge0} \int_{2^k}^{2^{k+1}} \frac{dx}{b(x)} \,\le\, \sum_{k\ge0} \frac{2^{k+1} - 2^k}{b(2^k)} \,=\, \frac14\sum_{k\ge0} \frac{1}{f(2^k)},
    \end{align*}
    and
    \begin{align*}
        \sum_{k\ge0} \int_{2^k}^{2^{k+1}} \frac{dx}{b(x)} \,\ge\, \sum_{k\ge0} \frac{2^{k+1} - 2^k}{b(2^{k+1})} \,=\, \frac{1}{8}\sum_{k\ge0} \frac{1}{f(2^{k+1})}.
    \end{align*}
   % Therefore
    %\begin{align*}
       % \frac{1}{2}\sum_{k\ge0} \frac{1}{f(2^{k+1})} \le \int_1^\infty \frac{dx}{b(x)} \le \sum_{k\ge0} \frac{1}{f(2^k)},
    %\end{align*}
    The lemma follows.    
\end{proof}

\begin{lemma} \label{lem:sab} Let $B$ be a Brownian motion, and for $a, b>0$ let 
\[ \sigma_{a,b}:= \inf\left\{t>0: B_t<bt-a\right\}.\]
We have 
\[\E\left[\exp\left(s\sigma_{a,b}\right)\right] = \exp\left(-ab\left[\sqrt{1-\frac{2s}{b^2}}-1\right]\right),\]
for $\frac{s}{b^2}<\frac12$.
\end{lemma}
\begin{proof} This follows from Exercise 3.14 in Chapter 2 of \cite{revu-yor}. Indeed, the hint given in the second part of the exercise shows that 
\[ \E\left[\exp\left(\frac{s}{b^2} \sigma_{ab,1}\right)\right] = \exp\left(-ab\left[\sqrt{1-\frac{2s}{b^2}}-1\right]\right), \]
as long as $\frac{s}{b^2}<\frac12$. It is easy to see by scaling properties that $\sigma_{ab,1} \stackrel{d}{=} \sigma_{a,b}$ and this completes the proof of the lemma. 
\end{proof}

\begin{proof}[Proof of Proposition \ref{prop:explosion}]  We will consider the evolution of the above SDE $X$ between levels $ 2^k, \, k \in \Z$. % where $a_k$ are such that
%\be \label{eq:ak:2} a_k = 2^k.\ee
%It is clear that $a_k$ is increasing with $k$.
 We will compare the evolution of $X$ to that of $Y^{(k)}$, which is the SDE
\[
dY^{(k)}_t = f(2^{k-1}) Y^{(k)}_t dt + Y^{(k)}_t dB_t \]
{\it starting} at $2^k$, whose solution is 
\[ Y^{(k)}_t = 2^k \exp\left( B_t + \left[f(2^{k-1})-\frac12\right]t\right).\]
It is clear that for $k$ large enough (so that $f(2^{k-1})-\frac12 >0$) the process $Y_t^{(k)}$ remains positive and goes to infinity as $t$ goes to infinity. %It is also clear from the comparison result for SDEs (Proposition 2.18 of Chapter 5 of \cite{kara-shre}) that $X$ starting at $2^k$, stays above $Y^{(k)}$ as long as $Y^{(k)}$ remains above $2^{k-1}$. 

Let us say that the process $X$ {\it backtracks at level $2^k$}, if starting at $2^k$ it hits level $2^{k-1}$ before level $2^{k+1}$. Fix $K>0$. 
\bes \begin{split}
&\bP_{2^K} \left(X \text{ explodes to } \infty \text{ by time } \delta\right) \\
& \ge \bP_{2^K} \Big( \Big\{X \text{ explodes to } \infty \text{ by time } \delta \text{ successively hitting levels } 2^k, \, k \ge K +1 \\
&\hspace{8cm}  \text{ without backtracking at any level}\Big\} \Big)
%&\hspace{5cm} \bigcap \Big\{ X \text{ explodes to } \infty \text{ by time } \delta\Big\} \Big) \\
%& \ge  P_{2^K} \Big( \Big\{X \text{ successively hits levels } 2^k, \, k \ge K  \text{ without backtracking at any level}\Big\}\Big) 
%& \hspace{2cm} + P_{2^K} \Big( \Big\{ X \text{ explodes to } \infty \text{ by time } \delta\Big\} \Big) -1
\end{split}
\ees
Define stopping times $S_0, S_1, \cdots$ as follows. Let $S_0=0$ and 
\[ S_{i+1} = \inf\left\{t>S_{i}: X_t \notin [2^{-1} X_{S_i}, 2X_{S_i}]\right\}.\]
Define $\{0,1\}$ valued random variables $J_i,\, i\ge 1$ as follows: $J_i$ equal to $1$ if $X_{S_i}=2X_{S_{i-1}}$, and equal to $0$ otherwise. Therefore %If $S$ is the exit time of the interval $[2^{k-1}, 2^{k+1}]$, then let $Z_2$ equal to $1$ if $X$ hits $2X_S$ before $2^{-1}X_S$, and $0$ otherwise. Define $Z_3, Z_4 \cdots$ recursively. Therefore 
\begin{align*}
&\bP_{2^K} \left(X \text{ explodes to } \infty \text{ by time } \delta\right) \\
&\ge \bP_{2^K} \left(\left\{\sum_{i=0}^{\infty} (S_{i+1}-S_i)\le \delta\right\} \bigcap\Big\{ J_1=J_2=\cdots =1\Big\}\right).
\end{align*}
Let us make an observation. Suppose $X_S= 2^k$ for some random stopping time $S$, and consider the SDE $Y^{(k)}$ driven by the Brownian motion $B_{S+\cdot}-B_S$ starting at $2^k$. It is clear by \eqref{eq:f}, the strong Markov property and the comparison  result for SDEs (Proposition 2.18 of Chapter 5 of \cite{kara-shre}) that $X_{S+\cdot}$ stays above $Y^{(k)}$ as long as $Y^{(k)}$ remains above $2^{k-1}$. 

Let us define an auxillary process $Z$ and stopping times $0=\tilde S_0\le  \tilde S_1\le \cdots$ as follows. The process starts at $2^K$ so that $Z_0=2^K$. Let 
\[dZ_t = f\left(\frac{Z_{\tilde S_i}}{2}\right) Z_t dt + Z_t dB_{t}, \quad \tilde S_i \le t< \tilde S_{i+1}, \]
with 
\[ \tilde{S}_{i+1} = \inf\left\{t>S_i: Z_t \notin [2^{-1} Z_{\tilde S_i}, 2Z_{\tilde S_i} ]  \right\}.\]
Similarly to earlier define $\{0,1\}$ valued random variables $\tilde J_i,\, i\ge 1$ as follows: $\tilde J_i$ equal to $1$ if $Z_{\tilde S_i}=2Z_{\tilde S_{i-1}}$, and equal to $0$ otherwise. Our observation in the previous paragraph gives 
\begin{align}
\label{eq:Z}
\begin{split}
&\bP_{2^K} \left(X \text{ explodes to } \infty \text{ by time } \delta\right) \\
&\ge \bP_{2^K} \left(\left\{\sum_{i=0}^{\infty} (\tilde S_{i+1}-\tilde S_i)\le \delta\right\} \bigcap\Big\{ \tilde J_1=\tilde J_2=\cdots =1\Big\}\right).
\end{split}
\end{align}
It is to important to note here that the above inequality holds because of the event $\tilde J_i=1$ for all $i$, as the comparison of $X$ with $Y^{(k)}$ only holds as long as $Y^{(k)}$ remains above $2^{k-1}$.

To proceed we compute the probability $\bP(\tilde J_i=1 \text{ for all } i)$. Let $T_{k-1,k+1}$ be the exit time of $Y^{(k)}$ from the interval $[2^{k-1}, 2^{k+1}]$.
It follows from equation (5.61) of Chapter 5 of \cite{kara-shre} that 
\bes
%P_{a_{k}}\left[Y_{T_{k-1,k+1}}=2^{k-1}\right]&= \frac{p(2^{k+1}) -p(2^{k})}{p(2^{k+1}) -p(2^{k-1})} \\
\bP_{2^{k}}\left[Y^{(k)}_{T_{k-1,k+1}}=2^{k+1}\right]= \frac{p(2^{k}) -p(2^{k-1})}{p(2^{k+1}) -p(2^{k-1})},
\ees
where
\[ p(x):= \int_{2^{k-1}}^x \exp\left(-2 \int_{2^{k-1}}^{\xi}\frac{f(2^{k-1}) \zeta }{\zeta^2} d\zeta\right)\, d \xi\]
is the {\it scale function}. Therefore 
\bes \label{eq:jump:prob}
\begin{split}
\bP_{2^{k}}\left[Y^{(k)}_{T_{k-1,k+1}}=2^{k+1}\right]&= \frac{\int_{2^{k-1}}^{2^{k}}\exp\left(-2 f(2^{k-1}) \log\left(\frac{\xi }{2^{k-1}}\right)\right)d\xi}{\int_{2^{k-1}}^{2^{k+1}}\exp\left(-2 f(2^{k-1})\log\left(\frac{\xi }{2^{k-1}}\right)\right)d\xi} \\
&= \frac{1}{1+ \frac{\int_{2^{k}}^{2^{k+1}}\left(\frac{2^{k-1}}{\xi}\right)^{2 f(2^{k-1}) }d\xi}{\int_{2^{k-1}}^{2^{k}}\left(\frac{2^{k-1}}{\xi}\right)^{2 f(2^{k-1}) }d\xi}}\\
&= \frac{1}{1+2^{1-2f(2^{k-1})}}
\end{split}
\ees
By the strong Markov property, we have for large $K$
\be \label{eq:succ_right}\begin{split} \bP_{2^K}\left(\tilde J_i=1 \text{ for all } i\right)&=\prod_{k \ge K} \bP_{2^{k}}\left[Y^{(k)}_{T_{k-1,k+1}}=2^{k+1}\right] \\&=\prod_{k \ge K}\frac{1}{1+2^{1-2f(2^{k-1})}}   \\
& \ge \prod_{k \ge K} \left[1- 2^{1-2f(2^{k-1})}\right] \\
& \ge 1- \sum_{k \ge K}  2^{1-2f(2^{k-1})},
\end{split}\ee
where we have used $(1+x)^{-1} \ge 1-x$  for all $x \ge 0$, and $(1-x)(1-y) \ge 1-x-y$ whenever $0\le x, y \le 1$. We note here that 
\bes \begin{split}  \sum_{k \ge K}  2^{1-2f(2^{k-1})} & \lesssim \sum_{k\ge K}\frac{1}{f(2^{k-1})}.
%& \le \sum_{k\ge K}\frac{2^{k-1}}{b(2^{k-1})}  <\infty,
\end{split}
\ees
By choosing $K$ to be large we can make the term on the left of \eqref{eq:succ_right} to be arbitrarily close to $1$ by Lemma \ref{lem:osgood}.

Next we note
\bes
\begin{split}
T_{k-1, k+1} &= \inf\Big\{t>0: Y^{(k)}_t \text{ leaves } [2^{k-1}, 2^{k+1}]\Big\} \\
&\le \inf\Big\{t>0: Y^{(k)}_t \text{ leaves }(0, 2^{k+1}]\Big\} \\
&= \inf\left\{t>0:  B_t + \left[f(2^{k-1})-\frac12\right]t> \log 2\right\}, \end{split}\ees
as the process $Y^{(k)}$ does not hit $0$ in finite time. Consequently 
\be \label{eq:tau:T} 
\begin{split}
T_{k-1, k+1} &\le \inf\left\{ t>0: B_t >-\left[f(2^{k-1})-\frac12\right]t+\log 2\right\} \\
&\stackrel{d}{=} \inf\left\{ t>0: B_t <\left[f(2^{k-1})-\frac12\right]t-\log 2\right\} \\
&=: T_k.
\end{split}
\ee
From Lemma \ref{lem:sab} 
\bes \label{eq:exp:T} \E\exp(sT_k) =\exp\left(-\left[f(2^{k-1})-\frac12\right]\log 2\cdot \left(\sqrt{1-\frac{2s}{\left[f(2^{k-1})-\frac12\right]^2}} -1\right)\right),\ees
as long as $\frac{s}{\left[f(2^{k-1})-\frac12\right]^2} <\frac12$. 
%From \eqref{eq:succ_right}, the first term is bounded below by $ 1- \sum_{k \ge K}  2^{1-2f(2^{k-1})}$.
%Let us make an observation. Suppose $X_S= 2^k$ for some random stopping time $S$, and consider the SDE $Y^{(k)}$ driven by the Brownian motion $B_{S+\cdot}-B_S$ starting at $2^k$. It is clear by the strong Markov property and the comparison  result for SDEs (Proposition 2.18 of Chapter 5 of \cite{kara-shre}) that $X_{S+\cdot}$ stays above $Y^{(k)}$ as long as $Y^{(k)}$ remains above $2^{k-1}$.

Consider {\it independent} random variables $\tilde T_{K}, \tilde T_{K+1},\cdots $ such that $\tilde T_k \stackrel{d}{=} T_k$ for all $k$.  From the above
\begin{align}
\label{eq:T:tail}
\begin{split}
&\bP(\tilde T_{K} +\tilde T_{K+1}+\cdots > \delta) \\
& \le e^{-s \delta } \E \exp\left[s\left(\tilde T_{K} +\tilde T_{K+1}+\cdots\right)\right] \\
&\le e^{-s\delta}\exp\left(-\sum_{k \ge K}\left[f(2^{k-1})-\frac12\right]\log 2\cdot \left(\sqrt{1-\frac{2s}{\left[f(2^{k-1})-\frac12\right]^2}} -1\right)\right) \\
&\le  e^{-s\delta}\exp\left(2s\log 2\cdot\sum_{k\ge K} \frac{1}{[f(2^{k-1})-\frac12]} \right),
\end{split}
\end{align}
using the inequality $\sqrt{1-x} \ge 1- x$, valid for $0\le x\le 1$. Returning to \eqref{eq:Z}, we obtain from \eqref{eq:succ_right}, \eqref{eq:tau:T}, \eqref{eq:T:tail} and the strong Markov property
\bes \begin{split}& \bP_{2^K} \left(X \text{ explodes to } \infty \text{ by time } \delta\right)  \\
&\ge \bP(\tilde T_{K+1} +\tilde T_{K+1}+\cdots \le \delta)+\bP_{2^K}\left(\tilde J_i=1 \text{ for all } i\right)-1 \\
& = \bP_{2^K}\left(\tilde J_i=1 \text{ for all } i\right)-  \bP(\tilde T_{K+1} +\tilde T_{K+1}+\cdots > \delta) \\
&\ge 1- \sum_{k \ge K}  2^{1-2f(2^{k-1})} - e^{-s\delta}\exp\left(2s\log 2\cdot\sum_{k\ge K} \frac{1}{[f(2^{k-1})-\frac12]} \right).
%& \ge P_{2^K}\left(\tilde J_i=1 \text{ for all } i\right)- e^{-s \delta } E \exp\left[s\left(\tilde T_{K} +\tilde T_{K+1}+\cdots\right)\right].
\end{split}
\ees
%Using \eqref{eq:exp:T} 
%\bes\begin{split} &P_{2^K} \left(X \text{ explodes to } \infty \text{ by time } \delta\right)  \\
%& \ge 1- e^{-s\delta}\exp\left(-\sum_{k \ge K}\left[f(2^{k-1})-\frac12\right]\log 2\cdot \left(\sqrt{1-\frac{2s}{\left[f(2^{k-1})-\frac12\right]^2}} -1\right)\right)
%\end{split}
%\ees
%The inequality $\sqrt{1-x} \ge 1- x$, valid for $0\le x\le 1$, gives
%\[P_{2^K} \left(X \text{ explodes to } \infty \text{ by time } \delta\right) \ge 1- e^{-s\delta}\exp\left(2s\log 2\cdot\sum_{k\ge K} \frac{1}{[f(2^{k-1})-\frac12]} \right). \]
The proof of Proposition \ref{prop:explosion} is complete by a large choice of $K$ using Lemma \ref{lem:osgood}. 
\end{proof}

\section{The Alternating Process}\label{sec:AP}
Let us recall the truncated interacting SDEs given in \eqref{eq:UJ} 
\be\label{eq:UJ:2} dU_t^{(J)}(x) =\left(\mathscr{L} U_t^{(J)}\right)(x) dt + b\left(U_t^{(J)}(x) \wedge J\right) dt + U_t^{(J)}(x)\, dB_t(x). \ee
We start with a {\it bounded nonnegative} initial profile $U_0$. The aim of this section is to compare the above process at a particular site to the underlying SDE. To this end, we consider the {\it Alternating Process} $V_t(x) = V_t^{(n,J)}(x),\, x\in \Z$ (for brevity, we will often drop the superscripts) defined as follows. The initial profile for this process is {\it the same }: $V_0=U_0$. Till time $\frac{1}{n}-$, the process is governed by independent SDEs:
\[ dV_t(x) = b\left(V_t(x) \wedge J\right) dt +V_t(x) \,dB_t\left(x\right),\quad x\in  \Z,\quad 0\le t<\frac1n.\]
At time $\frac1n$, the effect of the operator $\mathscr L$ kicks in:
\[ V_{\frac1n}(x) =\sum_{y \in  \Z} G\left(\frac1n; y-x\right) V_{\frac{1}{n}-}(y).\]
%where $G$ is the heat kernel for $\mathscr L$ {\color{red} (perhaps define $G$ earlier using random walk)}. 
This process is repeated. That is from time $\frac1n$ to time $\frac{2}{n}-$ we run independent SDEs with initial profile $V_{\frac1n}$:
\bes\begin{split}dV_t(x) = b\left(V_t(x) \wedge J\right) dt +V_t(x)\, dB_t\left(x\right),\quad x\in  \Z,\quad \frac1n\le t<\frac2n,\\
\end{split}\ees
 and 
 \[ V_{\frac2n}(x) =\sum_{y \in  \Z} G\left(\frac1n; y-x\right) V_{\frac{2}{n}-}(y),\]
and so on up to time $1$. The key result in this section is Proposition \ref{prop:VtoU} where we show that at any given site $x\in\Z$ the process $V^{(n,J)}_{\cdot}(x)$ converges as a process to $U^{(J)}_{\cdot}(x)$ (see also \cite{gyon-kryl}). This will allow us to compare $U^{(J)}$ to the underlying SDE (at that site).  In the following proposition, $[ x]$ will denote the largest integer less than or equal to $x$.

\begin{remark}
    The Alternating process can be defined when general Lipschitz continuous drifts and noise coefficients are introduced in the equation \eqref{eq:UJ:2}, and results prior to Proposition \ref{prop:compare} also hold at that level of generality when suitably restated. However, we only present the statements that are relevant for our purposes.
\end{remark}

\begin{proposition} The process $V_t(x)=V^{(n, J)}_t(x),\; x\in Z, \; 0\le t\le 1$ has the following integral formulation.
\be \label{eq:V}
\begin{split}
V_t(x) &= \left(G_{\frac{[nt]}{n}}*V_0\right) (x) + \int_0^{\frac{[nt]}{n}} \sum_{y \in \Z} G\left(\frac{1}{n}\left[ n\left(\frac{[nt]}{n}-s\right)+1\right]; y-x \right) b\left(V_s(y)\wedge J\right)ds \\
& \qquad +\int_0^{\frac{[nt]}{n}}\sum_{y \in  \Z} G\left(\frac{1}{n}\left[n\left(\frac{[nt]}{n}-s\right)+1\right]; y-x \right)  V_s(y) dB_s\left(y\right)  \\
& \qquad+ \int_{\frac{[nt]}{n}}^t b\left(V_s(x)\wedge J\right) ds +\int_{\frac{[nt]}{n}}^t V_s(x) dB_s\left(x\right)
\end{split}
\ee
\end{proposition}
\begin{proof}
We will prove this by induction. For $0\le t<\frac1n, \, x\in  \Z$
\[ V_t(x) = V_0(x) +\int_0^t b\left(V_s(x)\wedge J\right) ds +\int_0^tV_s(x) dB_s\left(x\right).\]
At time $t=\frac1n,$
\begin{align*}
\begin{split}
&V_{\frac1n}(x) \\
&=\sum_{y \in \Z} G\left(\frac1n;\, y-x\right) V_{\frac1n-}(y) \\
& = \sum_{y \in \Z} G\left(\frac1n;\, y-x \right) \left[V_0(y)+ \int_0^{\frac1n} b\left(V_s(y)\wedge J\right) ds + \int_0^{\frac1n}V_s(y) dB_s\left(y\right)\right] \\
& = \left(G_{\frac{1}{n}}*V_0\right) (x)+ \int_0^{\frac1n} \sum_{y \in  \Z} G\left(\frac1n; y-x \right)b\left(V_s(y)\wedge J\right) ds  + \int_0^{\frac1n} \sum_{y \in  \Z} G\left(\frac1n; y-x \right)V_s(y) dB_s\left(y\right),
\end{split}
\end{align*}
which matches with \eqref{eq:V} with $t=\frac1n$. 

Now assume \eqref{eq:V} holds for $t\le \frac{k}{n}$. For $\frac{k}{n} < t<\frac{k+1}{n}$, 
\bes
\begin{split}
V_t(x) &= V_{\frac{k}{n}}(x) + \int_{\frac{k}{n}}^t b\left(V_s(x)\wedge J\right) ds +  \int_{\frac{k}{n}}^tV_s(x)  dB_s\left(x\right) \\
&=  \left(G_{\frac{k}{n}}*V_0\right) (x) + \int_0^{\frac{k}{n}} \sum_{y \in \Z} G\left(\frac{1}{n}\left[n\left(\frac{k}{n}-s\right)+1\right]; y-x \right) b\left(V_s(y)\wedge J\right)ds \\
& \qquad \qquad +\int_0^{\frac{k}{n}}\sum_{y \in  \Z} G\left(\frac{1}{n}\left[n\left(\frac{k}{n}-s\right)+1\right]; y-x \right)  V_s(y) dB_s\left(y\right)   \\
&\qquad +  \int_{\frac{k}{n}}^t b\left(V_s(x)\wedge J\right) ds +  \int_{\frac{k}{n}}^tV_s(x)  dB_s\left(x\right). 
\end{split}
\ees
This matches with \eqref{eq:V} since $[nt]=k$. Finally for $t=\frac{k+1}{n}$,
\bes
\begin{split}
V_{\frac{k+1}{n}}(x)&=\sum_{y \in \Z} G\left(\frac1n;\, y-x \right) V_{\frac{k+1}{n}-}(y) \\
%&=\sum_{i \in \epsilon\Z} G^{(\epsilon)}\left(\frac1n;\, x, i \right) \left[ V^{(n)}_{\frac1n}(i) + \int_{\frac1n}^{\frac2n} \left[V_s^{(n)}(i)\wedge M\right]^{1+\alpha} ds + \int_{\frac1n}^{\frac2n} \frac{1}{\sqrt \epsilon} V_s^{(n)}(i) dB_s\left(\frac{i}{\epsilon}\right)\right] \\
&= \left(G_{\frac{k+1}{n}}*V_0\right) (x) + \int_0^{\frac{k}{n}} \sum_{y \in \Z} G\left(\frac1n+\frac{1}{n}\left[n\left(\frac{k}{n}-s\right)+1\right]; y-x \right) b\left(V_s(y)\wedge J\right)ds \\
& \qquad \qquad +\int_0^{\frac{k}{n}}\sum_{y \in  \Z} G\left(\frac1n+\frac{1}{n}\left[n\left(\frac{k}{n}-s\right)+1\right]; y-x \right)  V_s(y) dB_s\left(y\right)   \\
&\qquad + \int_{\frac{k}{n}}^{\frac{k+1}{n}} \sum_{y \in \Z} G\left(\frac1n; y-x \right)b\left(V_s(y)\wedge J\right) ds
 + \int_{\frac{k}{n}}^{\frac{k+1}{n}} \sum_{y \in \Z} G\left(\frac1n; y-x \right)V_s(y) dB_s\left(y\right).
\end{split}
\ees
This also matches with \eqref{eq:V} since $\frac1n+\frac{1}{n}\left[n\left(\frac{k}{n}-s\right)+1\right]= \frac1n\left[n\left(\frac{k+1}{n}-s\right)+1\right]$. The proof of the proposition is complete by induction on $k$. 
\end{proof}

The following result is standard. See for example \cite{geor-jose-khos-shiu}, \cite{shig-shim}, \cite{foon-khos-col}.
\begin{lemma} \label{lem:U:mom} There exists a positive constant $C=C(\|U_0\|_{\infty}, J)$ such that \[ \sup_{\stackrel{t\le 1}{x\in \Z}} \E\left[\left\{U_t^{(J)}(x)\right\}^2\right] \le C.\]
\end{lemma}

We also have moment bounds on the alternating process $V_t(x)= V_t^{(n,J)}(x)$. Recall that we drop the superscripts for brevity. 
\begin{lemma}\label{lem:V:mom} There exists a constant $C_1=C_1(\|V_0\|_{\infty}, J)$, independent of $n$, such that 
\be \label{eq:V:mom} \sup_{\stackrel{t\le 1}{x\in \Z}} \E\left[\left\{V^{(n,J)}_t(x)\right\}^2\right] \le C_1.\ee 
\end{lemma}
\begin{proof} Consider instead the alternating process $\widetilde V= \widetilde V^{(n,J)}$ starting at the constant initial profile $\widetilde{V}_0\equiv \|V_0\|_{\infty}$. By comparison arguments \cite{muel}, it is enough to show \eqref{eq:V:mom} with $\widetilde V$ instead of $V$.  By the construction of the process $\widetilde{V}$, and stationarity in $x$,  it is clear that 
\[ \sup_{\stackrel{t\le 1}{ x\in  \Z}}\E \left[\widetilde{V}_t(x)^2\right] <\infty.\]
It is not immediate that the bound is independent of $n$.  To see this, let us use the integral representation (equation \eqref{eq:V}) for $\widetilde{V}$. For a universal constant $C$, we have that $\E\left[\widetilde{V}_t(x)^2\right]$ is at most
\begin{align*}
\begin{split}
&C\left\{ \|V_0\|_{\infty}^2 + \E\left[\left\{\int_0^{\frac{[nt]}{n}} \sum_{y \in  \Z} G\left(\frac{1}{n}\left[n\left(\frac{[nt]}{n}-s\right)+1\right]; y-x \right) b\left(\widetilde{V}_s(y)\wedge J\right) ds\right\}^2\right] \right. \\
&\qquad\qquad  +\int_0^{\frac{[nt]}{n}}\sum_{y \in \Z} \left\{G\left(\frac{1}{n}\left[n\left(\frac{[nt]}{n}-s\right)+1\right]; y-x \right)\right\}^2  \E\left[\widetilde{V}_s(y)^2\right] ds \\
&\qquad\qquad\qquad\qquad\qquad \left. + \E\left[\left\{\int_{\frac{[nt]}{n}}^t b\left(\widetilde{V}_s(x)\wedge J\right) ds\right\}^2\right] +\int_{\frac{[nt]}{n}}^t \E\left[\widetilde{V}_s(x)^2\right] ds \right\}.
\end{split}
\end{align*}
Using the fact that $\widetilde{V}_s(x)$ is stationary in $x$, and $\sum_y  G(t; y-x) ^2 \le \sum_yG(t; y-x)=1 $, we obtain the following for all $t\le 1$:
\[ \E\left[\widetilde{V}_t(x)^2\right] \le C(\|V_0\|_{\infty},J) + C(J) \int_0^t \E\left[\widetilde{V}_s(x)^2\right] ds.\]
Gronwall's inequality completes the proof. 
\end{proof}

We next show that the alternating process $V^{(n,J)}$ converges to $U^{(J)}$ in $L^2$. Before that we will need a lemma.
\begin{lemma}\label{lem:g:diff} For $0\le t<t+h\le 1$ we have 
\[\sum_{x \in \Z} \left|G(t+h; x) - G(t; x) \right|  \lesssim h. \]
\end{lemma}
\begin{proof} We have 
\bes \begin{split} &\sum_{x \in \Z} \left|G(t+h; x) - G(t; x) \right|  \\
& =  \sum_{x \in \Z}  \left| \bP\left(\mathcal X_{t+h} = x\right)-  \bP\left(\mathcal X_{t} = x\right)\right| \\
& = \sum_{x\in \Z} \left| \bP\left(\mathcal X_{t} \neq x,\; \mathcal X_{t+h} = x\right)-  \bP\left(\mathcal X_{t} = x,\;\mathcal X_{t+h} \neq x\right)\right| \\
&\le \sum_{x \in  \Z} \sum_{y \in \Z, y\neq x} \bP\left(\mathcal X_{t} =y,\; \mathcal X_{t+h} = x\right) + \sum_{x \in \Z} \sum_{y \in \Z, y\neq x}\bP\left(\mathcal X_{t} =x,\; \mathcal X_{t+h} = y\right) \\
&=\sum_{y \in \Z} \sum_{x\in \Z, x\neq y} \bP\left(\mathcal X_{t} =y,\; \mathcal X_{t+h} = x\right) + \sum_{x \in  \Z} \sum_{y \in \Z, y\neq x}\bP\left(\mathcal X_{t} =x,\; \mathcal X_{t+h} = y\right) \\
&= \sum_{y \in  \Z} \bP\left(\mathcal X_{t} =y\right) \cdot O\left(h\right) + \sum_{x \in \Z} \bP\left(\mathcal X_{t} =x\right) \cdot O\left(h\right) \\
&= O\left(h\right).
\end{split} \ees
The second last equality follows from general properties of continuous time random walks. Indeed the probability that the random walk jumps at least once in an interval $[t, t+r]$ is $O(r)$.
\end{proof}
\begin{proposition}  \label{prop:V-U}We have the following.
\be \label{eq:V-U}\sup_{t\le 1}\sup_{x\in  \Z}\E \left[\left\{V^{(n,J)}_t(x) - U^{(J)}_t(x)\right\}^2\right]\stackrel{n \to \infty}{\longrightarrow} 0.\ee
\end{proposition}
\begin{proof}
%We first show 
%\bes \sup_{t\le 1} \sup_{x\in \epsilon \Z}\E \left[\left\{V^{(n)}_t(x) - U_t(x)\right\}^2\right]\stackrel{n \to \infty}{\longrightarrow} 0.\ees
Recall that both $V= V^{(n,J)}$ and $U^{(J)}$ are starting with the same initial profile $V_0=U_0$. Using the integral representations of $V$ and $U^{(J)}$, and lemmas \ref{lem:U:mom} and \ref{lem:V:mom} we obtain 
\bes
\begin{split}
&\E\left[\left \lbrace U^{(J)}_t(x)- V_t(x) \right\rbrace^2\right]\\
& \lesssim \frac{C(\|U_0\|_{\infty},J)}{n} +\left[\left(G_{t}*U_0\right) (x)-\left(G_{\frac{[nt]}{n}}*U_0\right) (x)\right]^2 \\
&\qquad+\E\left[\bigg\lbrace \int_0^{\frac{[nt]}{n}} \sum_{y \in  \Z} G\left(\frac{1}{n}\left[n\left(\frac{[nt]}{n}-s\right)+1\right]; y-x \right) b\left(V_s(y)\wedge J\right)ds \right.\\
&\left.\hspace{6cm}-\int_0^{\frac{[nt]}{n}} \sum_{y \in  \Z} G\left(t-s; y-x \right) b\left(U^{(J)}_s(y)\wedge J\right) ds\bigg\rbrace^2\right] \\
&  \qquad +\E\left[\bigg\lbrace\int_0^{\frac{[nt]}{n}}\sum_{y \in  \Z} G\left(\frac{1}{n}\left[n\left(\frac{[nt]}{n}-s\right)+1\right]; y-x \right)  V_s(y) dB_s\left(y\right) \right.\\
&\left.\hspace{6cm} - \int_0^{\frac{[nt]}{n}}\sum_{y \in \Z} G\left(t-s; y-x \right)  U^{(J)}_s(y) dB_s\left(y\right) \bigg\rbrace^2\right],
\end{split}
\ees
where the constant $\frac{C(\|U_0\|_{\infty},J)}{n}$ comes from the last two terms of \eqref{eq:V}. Therefore 
\bes
\begin{split}
&\E\left[\left \lbrace U^{(J)}_t(x)- V_t(x) \right\rbrace^2\right]\\
&\lesssim \frac{C(\|U_0\|_{\infty},J)}{n} +\left[\left(G_{t}*U_0\right) (x)-\left(G_{\frac{[nt]}{n}}*U_0\right) (x)\right]^2 \\
&\hspace{2mm} + \E\left[\left\lbrace \int_0^{\frac{[nt]}{n}} \sum_{y \in \Z} G\left(\frac{1}{n}\left[n\left(\frac{[nt]}{n}-s\right)+1\right]; y-x \right)\left( b\left(V_s(y)\wedge J\right) -b\left(U^{(J)}_s(y)\wedge J\right) \right)ds\right\rbrace^2\right] \\
&\hspace{2mm} + \E\bigg[\bigg\lbrace \int_0^{\frac{[nt]}{n}} \sum_{ \in  \Z} \left(G\left(\frac{1}{n}\left[n\left(\frac{[nt]}{n}-s\right)+1\right]; y-x \right)-G\left(t-s; y-x \right)\right)\cdot b\left(U^{(J)}_s(y)\wedge J\right)  \bigg\rbrace^2\bigg]\\
&\hspace{2mm} + \int_0^{\frac{[nt]}{n}} \sum_{y \in  \Z} \bigg\lbrace G\left(\frac{1}{n}\left[n\left(\frac{[nt]}{n}-s\right)+1\right]; y-x \right)\bigg\rbrace^2\cdot \E\left[\left\lbrace U^{(J)}_s(y)-V_s(y)\right\rbrace^2\right]ds \\
&\hspace{2mm} + \int_0^{\frac{[nt]}{n}} \sum_{y \in \Z} \left(G\left(\frac{1}{n}\left[n\left(\frac{[nt]}{n}-s\right)+1\right]; y-x \right)-G\left(t-s; y-x \right)\right)^2\cdot\E\left[U^{(J)}_s(y)^2\right] ds.
\end{split}
\ees
Using $\sum_x G(\cdot; x)=1$, $G(\cdot; x)^2 \le G(\cdot; x)$ and an application of Cauchy-Schwarz inequality to the first expectation on the right gives
\bes
\begin{split}
&\sup_{x\in \epsilon \Z}\E\left[\left \lbrace U^{(J)}_t(x)- V_t(x) \right\rbrace^2\right] \\
&\lesssim_{(\|U_0\|_{\infty}, J)} \frac{1}{n}+ \left(\sum_{y\in \Z} \left|G\left(t; y-x \right) -G\left(\frac{[nt]}{n}; y-x \right)\right|\right)^2
\\
&\qquad +  \int_0^{\frac{[nt]}{n}}  \sup_{x\in \Z}\E\left[\left\lbrace U^{(J)}_s(x)-V_s(x)\right\rbrace^2\right]ds \\
&\qquad + \bigg[\int_0^{\frac{[nt]}{n}} \sum_{y \in \Z} \left|G\left(\frac{1}{n}\left[n\left(\frac{[nt]}{n}-s\right)+1\right]; y-x \right)-G\left(t-s; y-x \right) \right|ds\bigg]^2 \\
&\qquad + \int_0^{\frac{[nt]}{n}} \sum_{y \in  \Z} \left(G\left(\frac{1}{n}\left[n\left(\frac{[nt]}{n}-s\right)+1\right]; y-x \right)-G\left(t-s; y-x \right)\right)^2 ds.
\end{split}
\ees
%Using the expression
%\[ G(t;x,i) = e^{-t}\sum_{l=0}^{\infty} \frac{t^l \,\bP(\mathcal Z_l=x-i)}{l!},\]
%where $p_l(x,i)$ is the probability of a {\color{red} simple random walk} on $\Z$ to jump from $x$ to $i$ in $l$ steps, 
%we easily see that 
%\[ G\left(\frac{1}{n}\left[n\left(\frac{[nt]}{n}-s\right)+1\right]; i-x \right) \stackrel{n\to\infty}{\longrightarrow}G\left(t-s; i-x \right).\]
Lemma \ref{lem:g:diff} as well as Gronwall's inequality then gives 
\[ \sup_{t\le 1}\sup_{x\in\Z} \E\left[\left \lbrace U^{(J)}_t(x)- V_t(x) \right\rbrace^2\right] \stackrel{n\to\infty}{\longrightarrow} 0.\]
This completes the proof. 
\end{proof}

The main result of this section is the  process level convergence of $V_t^{(n,J)}(x)$ to $U_t^{(J)}(x)$ at any point $x\in\Z$. 
\begin{proposition}\label{prop:VtoU} Fix $J >0$. For any $x\in\Z$ the process $\{V^{(n, J)}_t(x),\, 0\le t\le 1\}$ converges in distribution to the process $\{ U^{(J)}_t(x),\, 0\le t\le 1\}$ in the space $D[0,1]$ of càdlàg functions in $[0,1]$. 
\end{proposition}
\begin{proof}
Fix an $x\in\Z$. Proposition \ref{prop:V-U} guarantees the finite dimensional convergence of $V_t(x)=V^{(n,J)}_t(x)$ to $U^{(J)}_t(x)$. What remains to be shown is the relative compactness of the sequence of processes $V^{(n,J)}$. To this end we verify statement (b) of Theorem 3.8.6 of \cite{ethi-kurt}. As a first step we verify condition (8.37) of Theorem 3.8.8 of \cite{ethi-kurt}.  Now
\be \label{eq:Vth-Vt} \begin{split}
V_{t+h}(x)-V_t(x)&= \left[\left(G_{\frac{[n(t+h)]}{n}}*V_0\right)(x)  - \left(G_{\frac{[nt]}{n}}*V_0\right)(x)\right] \\
&\quad +\int_0^{\frac{[nt]}{n}}\sum_{y\in \Z}\bigg[G\left(\frac1n\left[n\left(\frac{[n(t+h)]}{n}-s\right)+1\right];y-x\right) \\
&\hspace{3cm}- G\left(\frac1n\left[n\left(\frac{[nt]}{n}-s\right)+1\right];y-x\right)\bigg]\cdot b\left(V_s(y) \wedge J\right) ds \\
&\quad + \int_{\frac{[nt]}{n}}^{\frac{[n(t+h)]}{n}}\sum_{y\in \Z} G\left(\frac1n\left[n\left(\frac{[n(t+h)]}{n}-s\right)+1\right];y-x\right)\cdot b\left(V_s(y) \wedge J\right) ds \\
&\quad + \int_0^{\frac{[nt]}{n}}\sum_{y\in \Z} \bigg[G\left(\frac1n\left[n\left(\frac{[n(t+h)]}{n}-s\right)+1\right];y-x\right) \\
&\hspace{3cm}- G\left(\frac1n\left[n\left(\frac{[nt]}{n}-s\right)+1\right];y-x\right)\bigg]\cdot V_s(y) \,dB_s\left(y\right) \\
&\quad + \int_{\frac{[nt]}{n}}^{\frac{[n(t+h)]}{n}}\sum_{y\in \Z} G\left(\frac1n\left[n\left(\frac{[n(t+h)]}{n}-s\right)+1\right];y-x\right)\cdot V_s(y) \, dB_s\left(y\right) \\
&\quad +\left[\int_{\frac{[n(t+h)]}{n}}^{t+h} b\left(V_s(x)\wedge J\right) ds - \int_{\frac{[nt]}{n}}^t b\left(V_s(x)\wedge J\right) ds\right] \\
&\quad +\left[\int_{\frac{[n(t+h)]}{n}}^{t+h}V_s(x) \, dB_s(x) - \int_{\frac{[nt]}{n}}^t V_s(x)\,dB_s(x)\right]\\
&=: I+II+III+IV+V+VI+VII,
\end{split} \ee
where the roman numerals represent the corresponding terms on the right. We have a similar decomposition of $V_t(x)-V_{t-h}(x)$; simply replace the $t+h$ in \eqref{eq:Vth-Vt} with $t$ and the $t$ in the equation with $t-h$. This gives rise to 
\[ V_t(x)- V_{t-h}(x) = I^*+II^*+III^*+IV^*+V^*+VI^*+VII^*,\]
where the terms on the right hand side are similarly defined. We need to show that for the function $q(x):= |x|\wedge 1$ we have
\[ \E\left[\left\{q\left(V_{t+h}(x) -V_t(x)\right)\right\}^{m}\cdot\left\{q \left(V_t(x)-V_{t-h}(x)\right)\right\}^{m}\right] \lesssim  h^{\theta}\]
for some $\theta>1$ and some positive integer $m$. 

Consider first the case when $h<\frac{1}{2n}$. When $\frac{k}{n}\le t<\frac{k}{n}+ \frac{1}{2n}$ for some $k$, the terms $I, II, III, IV, V$ are all $0$, the term $VI$ is of order $h$ and the second moment of term $VII$ is of order $h$. Burkholders's inequality guarantees that the fourth moment of term $VII$ are of order $h^{2}$. As a consequence, when $\frac{k}{n}\le t<\frac{k}{n}+ \frac{1}{2n}$ with $m=4$
\begin{align*}
\begin{split}
\E\left[\left\{q\left(V_{t+h}(x) -V_t(x)\right)\right\}^{4}\cdot\left\{q \left(V_t(x)-V_{t-h}(x)\right)\right\}^{4}\right] &\le \E\left[\left\{V_{t+h}(x) -V_t(x)\right\}^{4}\right] \\
&\lesssim_{(J, \|U_0\|_{\infty})} h^2.
\end{split}
\end{align*}
When $\frac{k}{n}+ \frac{1}{2n}\le t< \frac{k+1}{n}$ the terms $I^*, II^*, III^*, IV^*, V^*$ are all $0$, the term $VI^*$ is of order $h$, and the second moment of term $VII^*$ is of order $h$. A similar argument shows that
\begin{align*}
\begin{split}
\E\left[\left\{q\left(V_{t+h}(x) -V_t(x)\right)\right\}^{4}\cdot\left\{q \left(V_t(x)-V_{t-h}(x)\right)\right\}^{4}\right] &\le \E\left[\left\{V_{t}(x) -V_{t-h}(x)\right\}^{4}\right] \\
&\lesssim_{(J, \|U_0\|_{\infty})} h^2.
\end{split}
\end{align*}
%{\color{red} (We still have to worry about term $I$).}
In the case $h\ge \frac{1}{2n}$ we simply bound  
\[\E\left[\left\{q\left(V_{t+h}(x) -V_t(x)\right)\right\}^{4}\cdot\left\{q \left(V_t(x)-V_{t-h}(x)\right)\right\}^{4}\right] \le \E\left[\left(V_{t+h}(x) -V_t(x)\right)^4\right].\] % + \E\left[\left(V_{t}(0) -V_{t-h}(0)\right)^4\right]. \]
First of all we bound term $I$. From the boundedness of the initial profile and Lemma \ref{lem:g:diff}
\be \label{eq:I}
\begin{split}
\left|\left(G_{\frac{[n(t+h)]}{n}}*V_0\right)(x)  - \left(G_{\frac{[nt]}{n}}*V_0\right)(x)\right|  &\lesssim_{\|U_0\|_{\infty}} \sum_{y \in \Z} \left|G\left(\frac{[n(t+h)]}{n}; y-x\right) - G\left(\frac{[nt]}{n}; y-x\right)\right|\\
&\lesssim_{\|U_0\|_{\infty}} h. 
\end{split}
\ee
The integrals in terms $III, V, VI, VII$ are over intervals of length at most $4h$. As a result the fourth moments of these terms are of order $h^2$. Let us consider terms $II$ and $IV$ %(and similarly terms $II^*, IV^*$).
\bes
\begin{split}
II  &\lesssim_{(J)} \int_0^{\frac{[nt]}{n}} \sum_{i\in  \Z} \bigg|G\left(\frac1n\left[n\left(\frac{[n(t+h)]}{n}-s\right)+1\right];y-x\right) \\ &\hspace{5cm}- G\left(\frac1n\left[n\left(\frac{[nt]}{n}-s\right)+1\right];y-x\right) \bigg| ds
\end{split}
\ees
Because we are working in the case that $h\ge \frac{1}{2n}$ we easily see that 
\[ \frac1n\left[n\left(\frac{[n(t+h)]}{n}-s\right)+1\right] - \frac1n\left[n\left(\frac{[nt]}{n}-s\right)+1\right] \le 4h,\]
and therefore from Lemma \ref{lem:g:diff}
\[ II \lesssim_{(J)} h.\]
Similarly the second moment of term $IV$ is
\bes
\begin{split}
&\lesssim_{(J, \|U\|_\infty)} \int_0^{\frac{[nt]}{n}} \sum_{y\in \Z} \bigg|G\left(\frac1n\left[n\left(\frac{[n(t+h)]}{n}-s\right)+1\right];y-x\right) \\
&\hspace{6cm} - G\left(\frac1n\left[n\left(\frac{[nt]}{n}-s\right)+1\right];y-x\right) \bigg|^2 ds \\
&\lesssim_{(J, \|U\|_\infty)} \int_0^{\frac{[nt]}{n}} \sum_{y\in  \Z} \bigg|G\left(\frac1n\left[n\left(\frac{[n(t+h)]}{n}-s\right)+1\right];y-x\right) \\
&\hspace{6cm} - G\left(\frac1n\left[n\left(\frac{[nt]}{n}-s\right)+1\right];y-x\right) \bigg| ds \\
&\lesssim_{(J, \|U\|_\infty)} h,
\end{split}
\ees
%where the constant in the inequality depends on both $J$ and $\|U_0\|_{\infty}$ and 
since the summands are at most $1$, hence the fourth moment of term $IV$ is of order $h^2$. The only remaining part is to verify condition (8.30) in statement (b) of Theorem 3.8.6 in \cite{ethi-kurt}. That is, we need a uniform (in $n$) bound on $\E\left[\left|V_h(x) - V_0(x)\right|^8\right]$ which decays to $0$ as $h\to 0$. This can be proven just as above, and we leave the details to the reader. %It is in fact simpler since terms $II$ and $IV$ of \eqref{eq:Vth-Vt} are $0$ when $t=0$. 
\end{proof}

\begin{corollary} \label{cor:UV} Fix time $\delta>0$ and $K, J > 0$. For any $x\in\Z$ we have 
\[ \bP\left(\sup_{t\le \delta} U_t^{(J)}(x) \ge K\right) \ge \limsup_{n\to \infty} \bP\left(\sup_{t\le \delta} V_t^{(n,J)}(x)\ge K\right).\]
\end{corollary}
\begin{proof} Let us first show that the set $A:=\{X: \sup_{t\le \delta} X(t) \ge K$ is a closed set in $D[0,\delta]$ (the space of right continuous functions on $[0,\delta]$ with left limits) with Skorohod topology.   Suppose $f_n$ is a sequence of functions in $A$ converging to $f$ in $D[0,\delta]$. To show $A$ is closed we need to show that $f$ is in $A$. We use Proposition 5.3 in Chapter 3 of \cite{ethi-kurt}. Let $\Lambda$ be the family of strictly increasing mappings $\lambda(t)$ which are Lipschitz continuous, with $\lambda(0)=0$ and $\lambda(\delta)=\delta $. Since $f_n$ converges to $f$ there is a sequence $\{\lambda_n\} \subset \Lambda$ the following hold:
\begin{itemize}
 \item $\lim_{n\rightarrow\infty}\sup_{0\le t\le \delta} |\lambda_n(t) - t| =0$,
 \item $\lim_{n\rightarrow\infty}\sup_{0\le t\le \delta} |f(\lambda_n(t)) - f_n(t)| =0$.
\end{itemize}
% Fix any $\varepsilon>0$, find $N_0$ very large so that for all $n\ge N_0$, second condition gives
    %\begin{align*}
    %    \sup_{0\le t\le 1}|f(\lambda_n(t)) - f_n(t)| <\varepsilon.
  %  \end{align*}
 %   Rewrite 2 in the following way
 %   \begin{align*}
 %       f(\lambda_n(t))>f_n(t)-\varepsilon,
%    \end{align*}
   % for all $t\le 1$ and $n\ge N_0$. Taking supremum over $t\le1$ also varies $\lambda_n(t)$ over $[0,1]$, so we get
Fix any $\epsilon>0$. For large $n$    
\begin{align*}
        \sup_{t\le \delta}f(t) =\sup_{t\le\delta}f(\lambda_n(t))\ge\sup_{t\le1}f_n(t)-\epsilon\ge K-\epsilon.
    \end{align*}
    Since $\epsilon>0$ was arbitrary, $\sup_{t\le\delta}f(t) \ge K$, and hence $f$ lies inside $A$.  This proves that $A$ is closed. The corollary now follows from Theorem 3.1 in Chapter 3 of \cite{ethi-kurt}.
\end{proof}

\subsection{Comparison of the Alternating Process with the underlying SDE.} 
Fix a site $x\in\Z$. Due to \eqref{eq:b:cond2} we can find large enough $n_0$ and $K_b$ such that for all $n\ge n_0$
\be \label{eq:kb}
e^{-\frac{1}{n}}b(u) - u \ge \frac{b(u)}{n_0},\, \text{ for all } u\ge K_b.
\ee
Recall that $B_t(x)$ is the Brownian motion at site $x$; see \eqref{eq:ISDE}. Consider the process 
\be \label{eq:z:sde}
%\begin{split}
%dX_t &= b(X_t\wedge J) \,dt +  X_t \,dB_t(0), \\
dZ_t  = \frac{b (Z_t\wedge J)}{n_0}\, dt +Z_t\, dB_t(x),
%\end{split}
\ee
starting at initial value $Z_0= M$ for some $M>0$, and driven by the Brownian motion $B_t(x)$, %First note that by comparison principle $Z_t\le X_t$ for all time. 
and consider the alternating process $V^{(n,J)}$ starting with the initial profile $M\mathbf{1}_{\{x\}}(\cdot)$. %$ which equals $M$ at $x$ and is $0$ everywhere else.

%Let $K_b$ be the smallest positive number such that 
%\be \label{eq:kb} \frac{b(x)}{e} -x \ge \frac{b(x)}{3},\, \text{ for all } x\ge K_b.\ee

\begin{proposition}\label{prop:compare} Fix $M>0$ and $x$ as above. Assume that $b(0)=0$ and the condition \eqref{eq:b:cond2} holds. Fix $n_0$ as above. Then for every $n\ge n_0$, the process $V_t^{(n, J)}(x)$ stays above the process $Z_t$ in \eqref{eq:z:sde} as long as $K_b < Z_t < J$.
\end{proposition}
\begin{proof} 
Consider the process $X_t$ with $X_0=M$ that satisfies the SDE
\begin{align*}
    d X_t = b(X_t \wedge J)\,dt + X_t\, dB_t(x), \qquad \frac{l}{n} \le t < \frac{l+1}{n};
\end{align*}
for every non-negative integer $l$, and at the right end points $\frac{l}{n}$ we let%the process $X_t$ starts afresh with the initial values
\begin{equation}\label{Jump:reduce:X}
    %\begin{split}
  %  X_0 &= M, \qquad \text{when $l = 0$}; \\
    X_{\frac{l}{n}} = e^{-\frac{1}{n}}X_{\frac{l}{n}-}.%, \qquad \text{ $l\ge1$}.
  %  \end{split}
\end{equation}
Note that before time $\frac{1}{n}$, the alternating process $V_t(x)$ equals $X_t$ almost surely. At times $\frac{l}{n}$,
\begin{align*}\label{Jump:reduce:V}
    V_{\frac{l}{n}}(x) & = \sum_{y\in\Z} G\left(\frac{1}{n}; y-x \right) V_{\frac{l}{n}-}(y) \\
    & \ge G\left(\frac{1}{n}; 0 \right) V_{\frac{l}{n}-}(x) \\
    & \ge e^{-\frac{1}{n}}V_{\frac{l}{n}-}(x).
\end{align*}
Therefore $V_{\frac{1}{n}}(x) \ge X_{\frac{1}{n}}$. Now assume that $V_t(x)\ge X_t$ for all times $t\le \frac{l}{n}$. Since $V_t(x)$ and $X_t$ satisfy the same SDE in $[\frac{l}{n}, \frac{l+1}{n})$ it follows $V_t(x)\ge X_t$ for times $t< \frac{l+1}{n}$, and then the above argument gives $V_{\frac{l+1}{n}}(x) \ge X_{\frac{l+1}{n}}$. Therefore $V_t(x)\ge X_t$ for all time.

Define a new process $Y_t := e^{-\left(\frac{nt - [ nt ]}{n}\right)}X_t$. Observe that $Y_t$ is a continuous process due to \eqref{Jump:reduce:X}. Moreover the SDE satisfied by $Y_t$ is
\begin{align*}
    dY_t = \left[e^{-\frac{nt - [nt ]}{n}}\,b\left(e^{\frac{nt - [ nt ]}{n}} Y_t \wedge J\right) - Y_t \right ]\, dt + Y_t \,dB_t(x).
\end{align*}
Note that $X_t \ge Y_t$ trivially. Now define another process $Q_t$ satisfying the SDE
\begin{align*}
    dQ_t = \left(\frac{b(Q_t \wedge J)}{e^{\frac1n}} - Q_t\right)\, dt + Q_t\, dB_t(x),
\end{align*}
starting at $M$. When $Q_t$ hits $0$ both the drift and noise terms in the SDE for $Q_t$ become $0$, and therefore $Q_t$ must remain non-negative (see the comparison principle for SDEs: Propn. 2.18 in Chapter 5 of \cite{kara-shre}). Because $\frac{nt - [ nt ]}{n}\le \frac{1}{n}$ for all $n\ge1$, the drift of $Q_t$ is at most the drift of $Y_t$ in light of non-negativity of $Q_t$. Therefore once again the comparison principle establishes that $Y_t\ge Q_t$ for all times. Finally we note that due to the comparison principle and \eqref{eq:kb} the process $Q_t$ remains above the process $Z_t$ given by \eqref{eq:z:sde} as long as $K_b<Z_t<J$. All these considerations put together prove the proposition.
\end{proof}

\section{Spatial growth of the driftless interacting SDEs}\label{sec:ISDE:no_drift}
%{\color{red} Here we used the notation $G(t; x)$ while earlier we used $G(t;x,y)$. make it consistent}
The purpose of this section is to show that for each fixed $t>0$ the interacting SDEs \eqref{eq:ISDE} $U_{t}(x)$ can take arbitrarily large values as $x$ varies. By Mueller's comparison principle \cite{muel} it is sufficient to show this for the interacting SDEs $\widetilde U$ without drift:
\be
\label{eq:ISDE:2} d\widetilde U_t(x) =\left(\mathscr{L} \widetilde U_t\right)(x)\, dt  + \widetilde U_t(x) \,dB_t(x) ,\quad x \in \Z, 
\ee
with bounded initial profile $\wu_0(x)$. %We will further assume that the initial profile $\wu_0$ for equation \eqref{eq:ISDE:2} is bounded and satisfies the following symmetry, monotonicity and decay conditions
%\begin{align}
%    \label{u0:decay}
%    \wu_0(0) > \wu_0(x) = \wu_0(-x) \text{ for all $|x|>0$, }
%    \lim_{|x|\to\infty} \wu_0(x) = 0, \text{ and }
%    \lim_{|x|\to\infty} \frac{\left|\log \wu_0(x)\right|}{\sqrt{\log |x|}} = 0.
%\end{align}
Using arguments similar to \cite{chen-khos-kim} and \cite{conu-jose-khos} we will prove that for initial profiles that satisfy \eqref{u0:decay} the solution to \eqref{eq:ISDE:2} is unbounded in space for each time $t>0$.
Note that $\widetilde U$ is in fact just $U^{(0)}$.

%The arguments in this section follow closely that of \cite{conu-jose-khos} and \cite{chen-khos-kim}.
We first obtain upper and lower bounds on the moments of  $\widetilde U_{t}(0)$. 

\begin{lemma}\label{lem:mom:Ut} Consider the solution $\wu$ of \eqref{eq:ISDE:2} with initial profile $\wu_0 \equiv M$. For each $t \ge 0$, $x\in\Z$ and $k\ge 2$, the $k$-th moments of $\wu$ satisfy the bounds
\[ M^k \exp \left(\frac{k(k-3)}{2} t\right) \le \E\left[\left|\widetilde{U}_t(x)\right|^k\right] \le M^k2^{k/2}\exp\left(4k^2t\right).\]
\end{lemma}
\begin{proof} Assume the Lemma is true first for $M=1$. Then for any $M>0$ the process $U'_t(x) = \frac{\widetilde U_t(x)}{M}$ solves equation \eqref{eq:ISDE:2} with constant initial profile $\wu_0(x) \equiv 1$. The uniqueness of solution of \eqref{eq:ISDE:2} and the bounds for $U'_t(x)$ implies the Lemma. 

Let $M=1$ for the remainder of the proof. We first obtain an upper bound. For $x\in\Z$
\bes
\begin{split}
\widetilde{U}_t(x) &= 1 + \sum_{y\in\Z}\int_0^{t} G\left(t-s;y-x\right) \widetilde{U}_s(y)\, dB_s(y).
\end{split}
\ees
Define $\mathcal M_t^{(k)}:= \left\|\widetilde{U}_t(x)\right\|_k^2$. Using Burkholder's inequality we obtain for each 
\bes
\begin{split}
\mathcal M_t^{(k)} &\le 2 + 8k \int_0^t \sum_{y\in \Z} \left[G\left(t-s;y-x\right)\right]^2 \mathcal M_s^{(k)}\, ds \\
& = 2 + 8k \int_0^t \bP \left(\mathcal X_{2(t-s)} =0\right)M_s^{(k)}\, ds \\
& \le 2+ 8k \int_0^t \mathcal M_s^{(k)}\, ds, 
\end{split}
\ees
where $\mathcal X$ is a continuous time random walk with generator $\mathscr{L}$. Using Gronwall's inequality we obtain
\[\E\left[\left|\widetilde{U}_t(x)\right|^k\right]\le 2^{k/2}\exp\left(4k^2t\right).\]
We next obtain a lower bound on the moments. To obtain the lower bound we use the identity
\bes
\E\left[\left|\widetilde{U}_t(x)\right|^k\right] = \E \exp\left[\sum_{1\le i< j\le k} \int_0^t \mathbf{1}_{\{0\}}\left( \mathcal X^{(i)}_{s}-\mathcal X^{(j)}_{s}\right) ds\right],
\ees
where $\mathcal X^{(i)}, \, 1\le i \le k$ are independent continuous time random walks with generator $\mathscr L$. This is a result of the Feynman-Kac formula and is given on page 45 of \cite{carm-molc}. Note that \cite{carm-molc} gives the formula for the moments of the Stratonovich solution of \eqref{eq:ISDE:2} which will include an additional factor of $\exp(kt/2)$. A lower bound on the above is clearly
\[ \exp\left(\frac{k(k-1)t}{2}\right)\cdot \bP\left(\text{none of the $k$ walks $\mathcal X^{(i)}$ jump up to time $t$}\right).\]
Therefore we have 
\[ \E\left[\left|\widetilde{U}_t(x)\right|^k\right]  \ge  \exp \left(\frac{k(k-1)t}{2} - kt\right) \]
This completes the proof of the lemma. 
\end{proof}

The following ideas are inspired from \cite{conu-jose-khos}. For $\beta>0$ fixed, define 
\[ \widetilde{U}_t^{(\beta)}(x) := \sum_{y\in\Z} G(t;y-x)\wu_0(y)  + \sum_{y\in \Z: |y- x|\le\sqrt{\beta t}}\int_0^{t} G\left(t-s; y-x\right) \widetilde{U}^{(\beta)}_s(y) \, dB_s(y).  \]
The following lemma can be obtained after reworking the proof of upper bound in Lemma \ref{lem:mom:Ut}, hence we skip the details.
\begin{lemma} \label{lem:Ub:mom} Let $k\ge 2$, $t\ge0$ and $\beta>0$. If $\|\wu_0\|_\infty \le M$ then
\[ \sup_{x\in\Z} \E\left[\left|\widetilde{U}^{ (\beta)}_t(x)\right|^k\right] \le M^k 2^{k/2}\exp\left(4k^2t\right).\]
\end{lemma}
We next bound the moments of the difference of $\widetilde U_t(x)$ and $\widetilde U_t^{(\beta)}(x)$. The following lemma shows that if $\beta$ is large then $\wu_t(x)$ is close to $\wu_t^{(\beta)}(x)$. 
\begin{lemma} \label{lem:U-Ub} Fix $0<t\le 1$, $k\ge2$, $\beta>0$ and suppose $\|\wu_0\|_\infty \le M$. Let $\mathcal{R}$ be a bound on the jumps of the random walk $\mathcal{X}$ with generator $\mathscr{L}$. Then
\[ \sup_{x\in\Z} \E\left[\left\vert \widetilde{U}_t(x)-\widetilde{U}_t^{(\beta)}(x)\right\vert^k\right] \le  (8k)^{k/2} \exp\left(8k^2t\right) M^k \exp\left(\frac{k\sqrt{\beta t}}{\mathcal{R}} \log\left(\frac{\mathcal{R}e\sqrt t}{\sqrt \beta}\right)\right).\]
\end{lemma}
\begin{proof}
Let us define only for the purpose of the current proof 
\[\widetilde{V}_t^{(\beta)}(x) = \sum_{y\in\Z} G(t; y-x)\wu_0(y) +\sum_{y\in \Z}\int_0^{t} G\left(t-s; y-x\right) \widetilde{U}^{(\beta)}_s(y)\, dB_s(y). \]
By Burkholder's inequality and Lemma \ref{lem:Ub:mom}
\bes
\begin{split}
& \left\|\widetilde{V}_t^{(\beta)}(x)-\widetilde{U}_t^{(\beta)}(x)\right\|_k^2 \\
 &\le 4k \int_0^t \sum_{y:| y-x|>\sqrt{\beta t}} \left[G\left(t-s;y-x\right)\right]^2 \left\|\widetilde{U}^{(\beta)}_s(y)\right\|_k^2 ds  \\
& \le  8M^2k \exp\left(8kt\right) \int_0^t \sum_{y:| y|>\sqrt{\beta t}}\left[G\left(t-s;y\right)\right]^2 \,ds \\
& \le 8M^2k \exp\left(8kt\right) \int_0^t \bigg[ \sum_{y:| y|>\sqrt{\beta t}}G\left(s;y\right)\bigg]^2\,ds \\
&\le  8M^2k \exp\left(8kt\right)\cdot \bP\left(\text{RW $\mathcal X$ jumps more than $\frac{\sqrt{\beta t}}{\mathcal R}$ times in time $t$}\right)^2,
\end{split}
\ees
where $\mathcal X$ is a random walk with generator $\mathscr L$ with range $\mathcal R$. Since the jumps of the random walk happen at the times of a Poisson process with rate $1$, we consider $Y \sim \text{Poisson} (t)$ and compute $\bP(Y>\frac{\sqrt{\beta t}}{\mathcal R})$. For $\lambda>0$
\bes
\begin{split}
\bP\left(Y >\frac{\sqrt{\beta t}}{\mathcal R}\right) & \le \E\left[e^{\lambda Y - \lambda \frac{\sqrt{\beta t}}{\mathcal R}}\right] \\
& = \exp\left(te^{\lambda} -t -\lambda \frac{\sqrt{\beta t}}{\mathcal R}\right)
\end{split} 
\ees 
The right hand side is minimized for $\lambda= \log \left(\frac{\sqrt \beta}{\mathcal R \sqrt t}\right)$, and with this choice
\[  \bP\left(Y >\frac{\sqrt{\beta t}}{\mathcal R}\right) \le\exp\left(\frac{\sqrt{\beta t}}{\mathcal R} - t -\frac{\sqrt{\beta t}}{\mathcal R} \log \left(\frac{\sqrt\beta}{\mathcal R\sqrt t}\right)\right).\]
 The conclusion is 
\[ \bP\left(Y >\frac{\sqrt{\beta t}}{\mathcal R}\right) \le \left(\frac{\mathcal R e \sqrt t}{\sqrt \beta}\right)^{\frac{\sqrt{\beta t}}{\mathcal R}}.\]
In summary we have obtained 
\[ \sup_{x\in\Z} \left\|\widetilde{V}_t^{(\beta)}(x)-\widetilde{U}_t^{(\beta)}(x)\right\|_k^2 \le  8M^2k \exp\left(8kt\right) \left(\frac{\mathcal Re\sqrt t}{\sqrt \beta}\right)^{\frac{2\sqrt{\beta t}}{\mathcal{R}}}.\]
Next we bound $\left\|\widetilde{V}_t^{(\beta)}(x)-\widetilde{U}_t(x)\right\|_k^2$. It can be easily seen that 
\bes
\begin{split}
\left\|\widetilde{V}_t^{(\beta)}(x)-\widetilde{U}_t(x)\right\|_k^2 & \le 4k \int_0^t \sum_{y\in\Z} \left[G\left(t-s; y-x\right)\right]^2\left\|\widetilde{U}^{(\beta)}_s(y)-\widetilde{U}_s(y)\right\|_k^2 ds 
\end{split}
\ees
Therefore 
\bes
\begin{split}
&\sup_{x\in\Z} \left\|\widetilde{U}_t(x)-\widetilde{U}_t^{(\beta)}(x)\right\|_k^2\\
 & \le 8M^2k \exp\left(8kt\right) \left(\frac{\mathcal Re\sqrt t}{\sqrt \beta}\right)^{\frac{2\sqrt{\beta t}}{\mathcal R}}+ 8k \int_0^t \sum_{y\in\Z} \left[G\left(t-s; y-x\right)\right]^2\left\|\widetilde{U}^{(\beta)}_s(y)-\widetilde{U}_s(y)\right\|_k^2 ds  \\
& \le 8M^2k \exp\left(8kt\right) \left(\frac{\mathcal Re\sqrt t}{\sqrt \beta}\right)^{\frac{2\sqrt{\beta t}}{\mathcal R}}+ 8k \int_0^t \sup_{y\in\Z} \left\|\widetilde{U}^{(\beta)}_s(y)-\widetilde{U}_s(y)\right\|_k^2 ds 
\end{split}
\ees
Gronwall's inequality therefore gives us 
\[ \sup_{x\in\Z} \left\|\widetilde{U}_t(x)-\widetilde{U}_t{(\beta)}(x)\right\|_k^2 \le  8M^2k \exp\left(16kt\right) \left(\frac{\mathcal Re\sqrt t}{\sqrt \beta}\right)^{\frac{2\sqrt{\beta t}}{\mathcal R}}.\]
This completes the proof of the lemma.
\end{proof}

Next define $\widetilde{U}^{(\beta, n)}_t(x)$ to be the $n$th Picard iterate of $\widetilde{U}^{ (\beta)}_t(x)$: that is $\widetilde{U}^{ (\beta, 0)}_t \sum_{y\in\Z} G(t; y-x)\wu_0(y)$, and 
\bes
\widetilde{U}_t^{(\beta, n+1)}(x) := \sum_{y\in\Z} G(t; y-x)\wu_0(y) + \sum_{y\in \Z: | y- x|\le\sqrt{\beta t}}\int_0^{t} G\left(t-s; y-x\right) \widetilde{U}^{(\beta, n)}_s(y) dB_s(y).
\ees

\begin{lemma}\label{lem:U-Ubn} Assume conditions in Lemma \ref{lem:U-Ub}. If we choose any integer $n\ge 8ekt$ then we have
\be \label{eq:U-Ub}
\begin{split}
\sup_{x\in \Z} \E\left[\left|\widetilde{U}_t^{(\beta, n)}(x)- \widetilde{U}_t(x)\right|^k\right] &\lesssim 2^k M^k \left(\frac{1}{2}\right)^{\frac{k(n+1)}{2}} \\
&\hspace{0.5cm}+ M^k(32k)^{k/2} \exp\left(8k^2t\right) \exp\left(\frac{k\sqrt{\beta t}}{\mathcal R} \log\left(\frac{\mathcal Re\sqrt t}{\sqrt \beta}\right)\right). \\
%&\le 2^k (16k \epsilon)^{k/2} \exp(8k^2\epsilon) \exp\left(-C_5(1+\sqrt\beta)\cdot \frac{k}{2}\right) + (2C_6)^k \left[\sum_{m=n}^\infty m^{\frac14}\left(\frac{4ke\epsilon}{m}\right)^{\frac{m+1}{2}}\right]^k
\end{split}
\ee
%In particular when  $16k\epsilon \exp\left(-C_5(1+\sqrt \beta)\right) <\frac{1}{16}$ and $n \ge 8ek\epsilon$
%\be\label{eq:ub-ubn}
%\sup_{x\in \epsilon \Z}\sup_{t\le \epsilon^2} \E\left[\left|\widetilde{U}_t^{(\epsilon),(\beta, n)}(x)- \widetilde{U}_t^{(\epsilon)}(x)\right|^k\right] \le \frac{1}{2^k}\exp(8k^2\epsilon) +\frac{(C_7)^k}{2^{\frac{k(n+1)}{2}}}
%\ee
%In particular when $\beta= C_8(k\epsilon)^2$ for some large $C_8$ and $n \ge 8ek\epsilon$
%\be\label{eq:ub-ubn:2}
%\sup_{x\in \epsilon \Z}\sup_{t\le \epsilon^2} \E\left[\left|\widetilde{U}_t^{(\epsilon),(\beta, n)}(x)- \widetilde{U}_t^{(\epsilon)}(x)\right|^k\right] \le \exp(-2k^2\epsilon) +\frac{(C_7)^k}{2^{\frac{k(n+1)}{2}}}
%\ee
%{\color{red} One can get  weaker bound if $16k\epsilon \exp(16k\epsilon) \exp(-C_5(1+\sqrt\beta)) \le \frac{1}{16}$}
\end{lemma}
\begin{proof} Fix $0<t\le 1$ and $\beta>0$. We have for every $n\ge1$ and $k\ge2$ that
\begin{flalign*}
&\left\|\widetilde{U}_t^{(\beta, n+1)}(x) - \widetilde{U}_t^{(\beta, n)}(x)\right\|_k^2 &\\
&\qquad \le 4k \int_0^t \sum_{y\in \Z: | y -x|\le \sqrt{\beta t} }\left[G\left(t-s; y-x\right)\right]^2\left\|\widetilde{U}^{(\beta, n)}_s(y)-\widetilde{U}^{ (\beta, n-1)}_s(y)\right\|_k^2 ds.
\end{flalign*}
Defining
\[ \mathcal{M}_{n,k}(t):= \sup_{x\in\Z} \left\|\widetilde{U}_t^{(\beta, n+1)}(x)-\widetilde{U}_t^{(\beta, n)}(x)\right\|_k^2 ,\]
we obtain 
\[ \mathcal{M}_{n,k}(t) \le 4k \int_0^t \sum_{y\in \Z: | y -x|\le \sqrt{\beta t} }\left[G\left(t-s; y-x\right)\right]^2  \mathcal{M}_{n-1, k}(s)\, ds.\]
Note that since $\wu_0$ is bounded by $M$ it follows that  
\bes 
 \mathcal{M}_{0,k}(t) \le 4kM^2\int_0^t \sum_{y\in \Z}  \left[G\left(t-s; y-x\right)\right]^2 \, ds \le 4kM^2t.
\ees
From this we deduce by earlier arguments that 
\bes
\begin{split} \mathcal{M}_{n,k}(t) &\le 4k\int_0^t\mathcal{M}_{n-1,k}(s) ds \\
& \le \frac{M^2(4kt)^{n+1}}{n!},
\end{split}
\ees
In particular 
\[\sup_{x\in\Z} \left\|\widetilde{U}_t^{(\beta, n+1)}(x)-\widetilde{U}_t^{(\beta, n)}(x)\right\|_k \le  \frac{M}{\sqrt{n!}} (4kt)^{\frac{n+1}{2}}. \]
From the above it is clear that $\widetilde{U}_t^{(\beta, n)}(x)$ converges in $L^k$ to $\widetilde{U}_t^{(\beta)}(x)$ as $n \to \infty$. Furthermore for any $x\in\Z$
\bes \begin{split}\left\|\widetilde{U}_t^{(\beta)}(x)-\widetilde{U}_t^{(\beta, n)}(x)\right\|_k & \le \sum_{m=n}^{\infty} \left\|\widetilde{U}_t^{(\beta, m+1)}(x)-\widetilde{U}_t^{(\beta, m)}(x)\right\|_k \\
&\le {M} \sum_{m=n}^{\infty}\frac{1}{\sqrt{m!}} (4kt)^{\frac{m+1}{2}}  \\
&\lesssim {M} \sum_{m=n}^{\infty}m^{\frac14}\left(\frac{4kte}{m}\right)^{\frac{m+1}{2}},
\end{split}
\ees
where the last step follows by Sterling's approximation. If $n$ is chosen large enough so that $8ekt\le n$ then we have 
\[ \sup_{x\in\Z} \left\|\widetilde{U}_t^{(\beta)}(x)-\widetilde{U}_t^{(\beta, n)}(x)\right\|_k \lesssim M \left(\frac{1}{2}\right)^{\frac{n+1}{2}}.\]
The lemma follows from the above and Lemma \ref{lem:U-Ub}.
\end{proof}

The following lemma is proved similarly to Lemma 4.4. in \cite{conu-jose-khos}.
\begin{lemma} \label{lem:indep} Fix $0<t\le 1, \, \beta>0$ and $n\ge 0$. Fix $x_1, x_2,\cdots \in\Z$ such that $|x_i -x_j|\ge 2n(\sqrt{\beta t} +1) $ for all $i \ne j$. Then $\{ \widetilde{U}^{ (\beta, n)}_t(x_i)\}_{i \in \Z}$  is a collection of independent random variables. 
\end{lemma}

Let us state here a well known Bernstein-Freedman type tail probability estimate for the random walk $\mathcal{X}$ (see Theorem 2.9.5 in \cite{vers}). We give a proof for completeness.
\begin{lemma} \label{tail:G}
Let  $\mathcal X_t=\sum_1^{N_t} \mathcal Z_i$ be the random walk as in the introduction. %, i.e. $\mathcal Z_i$ are i.i.d. integer valued, symmetric random variables with variance $\sigma^2$ and bounded by $\RR$. 
Then for any $t, r>0$ we have
    \begin{align*}
        \bP(|\mathcal X_t|>r) \le 2\exp\left(-\frac{r^2}{2\left(\sigma^2t + \frac{r\mathcal{R}}{3}\right)}\right).
    \end{align*}
\end{lemma}
\begin{proof}
    We will use the inequality
    \begin{align*}
        e^{\lambda x} \le 1 + \lambda x + \frac{\lambda^2 x^2}{2(1-\frac{|\lambda| \RR}{3})},
    \end{align*}
    valid for all $|\lambda|<\frac{3}{\RR}, |x|\le\RR$ (so that $|\lambda x| < 3$). Symmetry of $\mathcal Z_1$ gives $\E[\mathcal{Z}_1]=0$. With this fact and the above inequality we obtain
    \begin{align}\label{Exp:estimate}
        \E\left[e^{\lambda \mathcal Z_1}\right] \le 1 + \frac{\lambda^2\sigma^2}{2\left(1-\frac{|\lambda| \RR}{3}\right)}.
    \end{align}
    Chernoff bound tells us that for any $\lambda>0$,
    \begin{align*}
        \bP(\mathcal X_t>r) &\le \E\left[ e^{\lambda \mathcal X_t} \right] e^{-\lambda r} \\
        &\le \exp \left( t( \E\left[ e^{\lambda \mathcal Z_1} \right] - 1) -\lambda r\right).
    \end{align*}
    Using inequality \eqref{Exp:estimate} we can improve the preceding into
    \begin{align*}
        \bP(\mathcal X_t>r) \le \exp \left( \frac{t\lambda^2\sigma^2}{2(1-\frac{|\lambda| \RR}{3})} - \lambda r \right).
    \end{align*}
    Therefore after plugging in $\lambda = \frac{r}{\sigma^2t+\frac{r\RR}{3}}<\frac{3}{\RR}$ gives
    \[
        \bP(\mathcal X_t>r) \le \exp \left( -\frac{r^2}{2\left(\sigma^2t+\frac{r\RR}{3}\right)} \right).
    \]
    Union bound gives the assertion.
\end{proof}
We now collect several preparatory Lemmas as in \cite{chen-khos-kim}.
\begin{lemma} \label{CalcLem}
    Fix integers $r$ and $ p$. Let $\wu_0 \in \ell^\infty(\Z)$ be a sequence which is zero in the interval $[p-r, p+r]$. Then for every $t>0$,
    \[
        \sup_{|x-p| \le r/2} \left|\sum_{y\in\Z} G(t; y-x)\wu_0(y)\right| \le 2 \|\wu_0\|_\infty \exp\left(-\frac{r^2}{8(\sigma^2t + r\RR/3)}\right)
    \]
\end{lemma}
\begin{proof}
    From the assumption on $\wu_0$ it follows that
    \[
        \left|\sum_{y\in\Z} G(t; y-x)\wu_0(y)\right| \le \|\wu_0\|_\infty \sum_{y:|y-p|>r} G(t; y-x) = \|\wu_0\|_\infty\sum_{z:|z+x-p|>r} G(t; z)
    \]
    Clearly $p-\frac{r}{2} \le x \le p+\frac{r}{2}$, and either $z+x>p+r$ or $z+x<p-r$ implies $|z|>r/2$. Applying the inequality in Lemma \ref{tail:G} with $r/2$ in place of $r$ proves the lemma.
\end{proof}
The following ideas are motivated by \cite{chendalang}. Define a new (commuting) binary operation `$\,\star\,$' on space-time indexed functions $f(t, x)$ and $g(t, x)$ for $t>0,\,x\in\Z$ as
\[
        (f \star g)(t, x) := \int_0^t \sum_{y\in\Z} f(t-s, x-y) g(s, y).
\]
Also define 
\begin{align*}%\label{eq:def:g}
\begin{split}
    g_0(t, x) &:=G^2(t; x), \\
    g_i(t, x) &:= (\underbrace{g_0\star g_0 \star \dots \star g_0}_{i+1})(t, x). %\qquad (i+1 \text{-fold convolution}).
\end{split}
\end{align*}
Let $h(t, x) = \sum_{i=0}^\infty g_i(t, x)$. It can be proven by induction that
\begin{align}\label{eq:g_i:bound}
    g_i(t, x) \le \frac{t^{i}}{i!}G^2(t;x).
\end{align}
Indeed,
\begin{align*}
    g_{i+1}(t, x) &= (\overbrace{g_0\star g_0 \star \dots \star g_{0} \star g_{0}}^{i+2})(t, x) \\%\qquad\qquad(\text{$i+2$ fold sum}) \\
    &\le \int_0^t \sum_{y\in\Z} \frac{(t-s)^i}{i!} G^2(t-s; x-y) G^2(s; y) ds \\
    &\le \int_0^t \frac{(t-s)^i}{i!} \left [\sum_{y\in\Z} G(t-s; x-y) G(s; y) \right]^2 ds \\
    &= \frac{t^{i+1}}{(i+1)!}G^2(t; x),
\end{align*}
by non-negativity of $G(t; x)$ and the Chapman-Kolmogorov equation. Therefore from \eqref{eq:g_i:bound} it follows that
\begin{align} \label{eq:h:bound}
    h(t, x) \le e^t G^2(t; x).
\end{align}
The next result puts these definitions in perspective. %We utilize the methods in \cite{chendalang}.
\begin{proposition} \label{NonGronwall}
    Fix $t, x$ and $\wu_0\in \ell^\infty(\Z)$. Let $q(t, x) := \sum_{y\in\Z} G(t; y-x)\wu_0(y)$. %Let $z(t, x)$ be any random field indexed by $(t, x)\in\R_+\times\Z$ which is at most {\color{red}exponentially growing in time $t$ and bounded in $x$}. If we have the inequality
    %\[
    %    z(t, x)\le q^2(t, x) + (z\star G)(t, x),
    %\]
    %for each $t$ and $x$ 
    Then
    \[
        \E\left[\left| \wu_t(x) \right|^2\right] \le q^2(t, x) + (q^2\star h)(t, x).
    \]
\end{proposition}
\begin{proof}
    Define the iterates $z^{(0)}(t, x)=q(t, x)$, and %From Lemma \ref{lem:g:diff} it follows that $z^{(0)}(t, x)$ is continuous in time due to boundedness of $f$. 
    %Let $z^{(1)}(t, x) = q(t, x) + \int_0^t \sum_{y\in\Z} G(t-s; x-y) z^{(0)}(s, y) dB_s(y) $ and apply Ito isometry to obtain
    %\begin{align*}
    %    \E\left[\left|z^{(1)}(t, x)\right|^2\right] &= q^2(t, x) + (q^2 \star g_0)(t, x) \\
    %    &\le q^2(t, x) + (q^2 \star h)(t, x),
    %\end{align*}
    %with $g_0$ and $h$ defined above. 
    for all $n\ge 1$
    \begin{align}\label{interates:z^k}
        z^{(n)}(t, x) = q(t, x) + \int_0^t \sum_{y\in\Z} G(t-s; x-y) z^{(n-1)}(s, y) dB_s(y).
    \end{align}
    Assume that
    \begin{align}\label{assume1z}
        \E\left[\left|z^{(n)}(t, x)\right|^2\right] &= q^2(t, x) + \sum_{i=0}^{n-1} (q^2 \star g_i)(t, x).
    \end{align}
    Hence the definition of $h$ gives
    \begin{align}\label{assume2z}
        \E\left[\left|z^{(n)}(t, x)\right|^2\right] \le q^2(t, x) + (q^2 \star h)(t, x).
    \end{align}
    %Observe that
    %\begin{align*}
    %    \E\left[\left|\int_0^t \sum_{y\in\Z} G(t-s; x-y) z^{(n)}(s, y) dB_s(y) \right|^2\right] &= \left( g_0 \star \E\left[\left|z^{(n)}\right|^2\right]\right)(t, x) \\
    %    &\le g_0 \star \left\{q^2(t, x) + (q^2 \star h)(t, x) \right\} \\
    %    &= (q^2 \star h)(t, x),
    %\end{align*}
    %since $h\star g_0 = h - g_0$. Therefore from equation \eqref{interates:z^k} and Ito isometry
    %\begin{align*}
    %    \E\left[\left|z^{(n+1)}(t, x)\right|^2\right] &\le q^2(t, x) + (q^2 \star h)(t, x).
    %\end{align*}
    From equations \eqref{interates:z^k} and \eqref{assume1z} we get by using Ito's isometry
    \begin{align*}
        \E\left[\left|z^{(n+1)}(t, x)\right|^2\right] &= q^2(t, x) + \left( g_0 \star \E\left[\left|z^{(n)}\right|^2\right]\right)(t, x) \\
        &= q^2(t, x) + \left( g_0 \star \left\{ q^2 + \sum_{i=0}^{n-1} (q^2 \star g_i) \right\}\right)(t, x) \\
        &=q^2(t, x) + \sum_{i=0}^{n} (q^2 \star g_i)(t, x),
    \end{align*}
    verifying \eqref{assume1z} and as a result \eqref{assume2z} for all $n$ by induction. Define 
    \begin{align*}
        \mathcal{M}^{(n)}(t, x) := \E \left[ \left| z^{(n+1)}(t, x) - z^{(n)}(t, x) \right|^2\right].
    \end{align*}
    From \eqref{interates:z^k} we find
    \begin{align*}
        \mathcal{M}^{(n)}(t, x) = \mathcal{M}^{(n-1)} \star g_0
    \end{align*}
    and
    \begin{align*}
        \mathcal{M}^{(0)}(t, x) = (q^2 \star g_0)(t, x).
    \end{align*}
    It thus follows from \eqref{eq:g_i:bound} that
    \begin{align*}
        \sum_{n=0}^\infty \left|\mathcal{M}^{(n)}(t, x)\right|^{\frac12} = \sum_{n=0}^\infty \left| (q^2 \star g_n)(t, x) \right|^\frac12 \le \sqrt{(q^2 \star g_0)(t,x)}\sum_{n=0}^\infty \sqrt{\frac{t^n}{n!}} <\infty.
    \end{align*}
    Therefore $z^{(n)}(t, x)$ is a Cauchy sequence in $L^2$. Let $z(t, x)$ be the limit. We need to verify that $z(t, x)$ coincides with our solution $\wu_t(x)$. %Once again by Ito isometry, we obtain
    %\begin{align*}
    %    \E \left[\left| z(t, x) - z^{(n+1)}(t, x) \right|^2\right] = \int_0^t \sum_{y\in\Z} G^2(t-s; y-x) \E \left[\left| z(t, x) - z^{(n)}(t, x) \right|^2\right] ds
    %\end{align*}
    %Taking $n \to \infty$ in \eqref{assume2z} and using Dominated convergence theorem gives
    Observe again by Ito isometry
    \begin{align*}
        \int_0^t \sum_{y\in\Z} G(t-s; x-y) z^{(n)}(s, y) dB_s(y) \;\stackrel{L^2}{\longrightarrow} \;\int_0^t \sum_{y\in\Z} G(t-s; x-y) z(s, y) dB_s(y).
    \end{align*}
    Consequently taking $n \to \infty$ in \eqref{interates:z^k} establishes the equality
    \begin{align*}
        z(t, x) = q(t, x) + \int_0^t \sum_{y\in\Z} G(t-s; x-y) z(s, y) dB_s(y),
    \end{align*}
    whence $z(t, x)$ is a solution of \eqref{eq:ISDE:2}. After taking limit $n \to \infty$ in equation \eqref{assume2z} we obtain
    \begin{align*}
        \E\left[\left| z(t, x) \right|^2\right] \le q^2(t, x) + (q^2\star h)(t, x),
    \end{align*}
    so the bound $q^2(t, x) \le \|\wu_0\|_\infty^2$ and inequality \eqref{eq:h:bound} shows
    \begin{align*}
        \sup_{t\le T,\, x\in\Z} \E\left[\left| z(t, x) \right|^2\right] \le 2\|\wu_0\|_\infty^2 e^{T} <\infty.
    \end{align*}
    As the equation \eqref{eq:ISDE:2} has the unique $L^2$ bounded solution $\wu_t(x)$ \cite[Theorem 2.1]{geor-jose-khos-shiu}, it is a modification of $z(t, x)$. The proof of proposition is now complete.
\end{proof}
%We skip the proof since it follows from modifying arguments as in \cite{chendalang}. 
We will also need the following.
\begin{lemma} \label{RobustISDE}
    Let $p\in\Z$ and $r>0$ be fixed. Let $\wu$ % and $\wv$ be two 
    solution to \eqref{eq:ISDE:2} with initial profile $\wu_0$  and assume that $\wu_0(x)\equiv 0$ for every $x\in[p-r, p+r]$. Then
    \[
        \sup_{x:|x-p|< r/4} \E\left[ \left|\wu_t(x)\right|^2 \right] \le 8\|\wu_0\|^2_\infty \exp\left(-\frac{r^2}{32\left(\sigma^2t + \frac{r\RR}{3} \right)} + t\right).
    \]
\end{lemma}
\begin{proof}
    Let $q(t, x) = \sum_{y\in\Z} G(t; y-x)\wu_0(y)$. Ito isometry gives
    \[
        \E\left[ \left|\wu_t(x)\right|^2 \right] \le q^2(t, x) + \int_0^t \sum_{y\in\Z} G^2(t-s; y-x) \E\left[ \wu^2_s(y) \right] ds.
    \]
    Therefore we can use symmetry of $G(t; x)$ and Proposition \ref{NonGronwall} and \eqref{eq:h:bound} to obtain
    \begin{align} \label{Lem2eq}
    \begin{split}
        \E\left[ \left|\wu_t(x)\right|^2 \right] &\le q^2(t, x) + (q^2\star h)(t, x) \\
        &\le q^2(t, x) + \int_0^t e^s \sum_{y\in \Z} q^2(t-s, x-y)G^2(s; y) ds.
    \end{split}
    \end{align}
    We split the summation inside the integral in the preceding as
    \begin{align*}
        &\sum_{y\in \Z} q^2(t-s, x-y)G^2(s; y) \\
        &= \sum_{|x-y-p|\le r/2} q^2(t-s, x-y)G^2(s; y) + \sum_{|x-y-p|> r/2} q^2(t-s, x-y)G^2(s; y) \\
    \end{align*}
    The first term can be simplified using Lemma \ref{CalcLem} as
    \begin{align*}
        \sum_{y:|x-y-p|\le r/2} q^2(t-s, x-y)G^2(s; y) &= \sum_{y:|x-y-p|\le r/2} q^2(t-s, x-y)G^2(s; y) \\
        &\le \sup_{z:|z-p|\le r/2} q^2(t-s, z) \left[ \sum_{|x-y-p|\le r/2} G^2(s; y) \right] \\
        &\le 4 \|\wu_0\|^2_\infty \exp\left(-\frac{r^2}{4(\sigma^2|t-s| + r\RR/3)}\right) \left[\sum_{y\in \Z} G(s; y) \right]^2.
    \end{align*}
    The second term can also be handled in a similar manner to get
    \begin{align*}
        \sum_{|x-y-p| > r/2} q^2(t-s, x-y) G^2(s; y) &\le \sum_{|x-y-p| > r/2} \left[ \sum_{z\in\Z} G(t-s; z-x+y)\wu_0(z) \right]^2 G^2(s; y) \\
        &\le \|\wu_0\|^2_\infty \left[ \sum_{z\in\Z} G(t-s; z-x+y) \right] \sum_{|x-y-p| > r/2} G^2(s; y)
    \end{align*}
    Therefore taking supremum over $|x-p|\le r/4$ after putting together the above calculations gives
    \begin{align} \label{Lem3eq1}
    \begin{split}
    &\sup_{x:|x-p|\le r/4} \sum_{y\in \Z} q^2(t-s, y-x)G^2(s; y) \\
    &\le 4 \|\wu_0\|^2_\infty \exp\left(-\frac{r^2}{4\left(\sigma^2t + \frac{r\RR}{3}\right)}\right) + 2\|\wu_0\|^2_\infty \left[\sum_{|y|\ge r/4} G(s; y) \right]^2 \\
    &\le 6 \|\wu_0\|^2_\infty \exp\left(-\frac{r^2}{32\left(\sigma^2t + \frac{r\RR}{3} \right)}\right),
    \end{split}
    \end{align}
    where we used Lemma \ref{tail:G}. Combining \eqref{Lem2eq} and \eqref{Lem3eq1} and another appeal to Lemma \ref{tail:G} finishes the proof.
\end{proof}
The next Lemma is of essential for our purposes and gives a control on the tail probabilities $\wu_t(x)$.
\begin{lemma} \label{log:wu:ratio}
    Let $\wu $ be the solution to \eqref{eq:ISDE:2} with initial profile $\wu_0$ satisfying the assumption \eqref{u0:decay}. For any fixed $\epsilon>0$ we have
    \[
        \lim_{|x|\to\infty} \frac{\log \bP(\wu_t(x) > \epsilon)}{\log |x|} =0. %\ge -\frac{Lc}{\sqrt{t}}.
    \]
\end{lemma}
\begin{proof}
    Since logarithm of probabilities is negative, the claim is reduced to proving that the limit inferior of the quantity in the claim is lower bounded by $0$.

    Let $p>0$ be an integer. Let $\wu^{(1), (p)}_0(x)=\wu_0(|x|\vee 4p)$, and $\wu^{(2), (p)}_0(x)=\wu_0(4p)$. Let $\wu^{(i), (p)}_t(x)$ be the solutions to the equation \eqref{eq:ISDE:2} with the initial profiles $\wu^{(i), (p)}_0(x)$, respectively, for $i=1, 2$. From assumptions \eqref{u0:decay} on $\wu_0$ and comparison principle we have the following for all time
    \[
    \wu^{(1), (p)}_t(x) \le \wu_t(x),
    \]
    and
    \[
    \wu^{(1), (p)}_t(x) \le \wu^{(2), (p)}_t(x).
    \]
    Lemma \ref{lem:mom:Ut} %, uniqueness (see Theorem 2.1 in \cite{geor-jose-khos-shiu}) and linearity of \eqref{eq:ISDE} 
    assures us that for each $k\ge 2$ and $t>0$ the following moment estimate holds for any $x\in\Z$:
    \[
    \wu_0^k(4p)\exp \left(\frac{k(k-3)}{2} t\right) \le \E\left[\left|\widetilde{U}^{(2), (p)}_t(x)\right|^k\right] \le \wu_0^k(4p) 2^{k/2}\exp\left(4k^2t\right).
    \]
    Therefore appealing to Paley-Zygmund inequality
    \begin{align*}
        \bP\left( \wu^{(2), (p)}_t(x) \ge \frac{\wu_0(4p)}{2}\exp\left[{\frac{k-3}{2}t}\right] \right) &\ge \bP \left( \wu^{(2), (p)}_t(x) \ge \frac{1}{2}\left\|\wu^{(2), (p)}_t(x)\right\|_{k} \right) \\
        &\ge \frac{1}{4} \frac{\E\left[ \left|\wu^{(2), (p)}_t(x) \right|^k \right]^2}{\E\left[ \left|\wu^{(2), (p)}_t(x) \right|^{2k} \right]} \\
        &\ge \frac{1}{4} \frac{\exp \left(\left[k(k-3) -16k^2\right]t\right) }{2^{k}} \\
        &\ge \frac{\exp \left(-17k^2t\right) }{2^{k+2}}.%,
    \end{align*}
    %for any $p\ge 2$. 
    Fix any $\epsilon>0$. Let $p$ and $k$ be such that
    \[
        \frac{\wu_0(4p)}{2}\exp\left[{\frac{k-3}{2}t}\right] = 2\epsilon,
    \]
    and note that we can ensure that $k \ge 2$ by making $p$ sufficiently large since $\wu(4p)$ goes to $0$ as $|p|$ increases to $\infty$. Thus with this choice of $p, k$ and using the inequality $(a+b)^2 \le 2a^2+2b^2$ we obtain
    \begin{align} \label{z>2ep}
    \inf_{x\in\Z} \bP\left( \wu^{(2), (p)}_t(x) \ge 2\epsilon \right) \ge A_1(t, \epsilon) \exp\left[ - \frac{136}{t} \log^2 \wu_0(4p) + \frac{\log 4}{t}\log \wu_0(4p) \right],
    \end{align}
    %using the value of $k$, the inequality $(a+b)^2\le 2a^2+b^2$, 
    where
    \[
    A_1(t, \epsilon)=\exp \left[-34t\left(\frac{2\log (4\epsilon)}{t}+3 \right)^2 - \frac{\log 4}{t}\log (4\epsilon)-\log 32\right].
    \]
    Also Markov's inequality, linearity of the equation \eqref{eq:ISDE:2}, and the Lemma \ref{RobustISDE} with $r=p/2$ gives
    \begin{align} \label{difference<ep}
    \begin{split}
        \sup_{x:|x-p|<p/8} \bP \left( | \wu^{(2), (p)}_t(x) - \wu^{(1), (p)}_t(x) | > \epsilon \right) &\le \sup_{x:|x-p|<p/8} \frac{\E \left[ | \wu^{(2), (p)}_t(x) - \wu^{(1), (p)}_t(x) |^2 \right]}{\epsilon^2} \\
        &\le 8\|\wu_0\|^2_\infty e^{t-2\log\epsilon}\exp\left(-\frac{p^2}{128(\sigma^2t + p\RR/6)} \right).
    \end{split}
    \end{align}
    Putting together \eqref{z>2ep} and \eqref{difference<ep} gives
    \begin{align*}
        \inf_{x:|x-p|<p/8} \bP\left( \wu_t(x) \ge 2\epsilon \right) &\ge \inf_{x:|x-p|<p/8} \bP\left( \wu^{(1), (p)}_t(x) \ge 2\epsilon \right) \\
        &\ge \inf_{x\in\Z} \bP\left( \wu^{(2), (p)}_t(x) \ge 3\epsilon \right) \\
        &- \sup_{x:|x-p|<p/8} \bP \left( | \wu^{(2), (p)}_t(x) - \wu^{(1), (p)}_t(x) | > \epsilon \right) \\
        &\ge A_1(t, \epsilon) \exp\left( - \frac{136}{t} \log^2 \wu_0(4p) + \frac{\log 4}{t}\log \wu_0(4p) \right) \\
        &\qquad\qquad\qquad - A_2(t, \epsilon)\exp\left(-\frac{p^2}{128(\sigma^2t + p\RR/6)} \right),
    \end{align*}
    where $A_2(t, \epsilon) = 8\|\wu_0\|^2_\infty e^{t-2\log\epsilon}$. 
    
    Fix any $0<\xi<1$. From assumptions \eqref{u0:decay} we can ensure that for large enough $p$ that 
    \begin{align*}
        %0< |\log \wu_0(4p)| < 1 \text{ and } 
        0>\log \wu_0(4p) \ge -\xi \sqrt{\log 4p}.
    \end{align*}
    Thus we obtain
    \begin{align*}
        \inf_{x:|x-p|<p/8} \bP\left( \wu_t(x) \ge 2\epsilon \right) &\ge A_1(t, \epsilon) \exp\left( - \frac{136\xi^2}{t} \log 4p - \frac{(\log 4)\xi}{t} \sqrt{\log 4p} \right) \\
        &\qquad\qquad- A_2(t, \epsilon)\exp\left(-\frac{p^2}{128(\sigma^2t + p\RR/6)} \right) \\
        &\ge A_1(t, \epsilon) \exp\left( - \frac{200\xi}{t} \log p \right) \\
        &\qquad\qquad- A_2(t, \epsilon)\exp\left(-\frac{p^2}{128(\sigma^2t + p\RR/6)} \right) \\
        &\ge \frac{A_1(t, \epsilon)}{2} \exp\left( - \frac{200\xi}{t} \log p \right)
    \end{align*}
    by choosing $p$ large enough. Therefore
    \begin{align*} 
        \liminf_{x \to \infty} \frac{\log \bP\left( \wu_t(x) \ge 2\epsilon \right)}{\log x} \ge - \frac{200\xi}{t},
    \end{align*}
    and the claim now follows as $0<\xi<1$ is arbitrary. % and $\wu_t(x)$ is symmetric.
\end{proof}

Now we can state the main result of this section.
\begin{proposition}\label{cor:sp:growth}
    Let $\wu$ be the solution of $\eqref{eq:ISDE:2}$ with non-negative, bounded initial profile $\wu_0$ that satisfies the decay condition \eqref{u0:decay}. Then
    \[
        \bP\left ( \sup_{x\in\Z} \wu_t(x) = +\infty \right) = 1.
    \]
\end{proposition}
\begin{proof}
     Let us fix a $K>0$ and $t>0$. Recall $\widetilde{U}^{ (\beta, 0)}_t = \sum_{y\in\Z} G(t; y-x)\wu_0(y)$, and 
    \bes
    \widetilde{U}_t^{(\beta, n+1)}(x):= \sum_{y\in\Z} G(t; y-x)\wu_0(y)
    +\sum_{y\in \Z: | y- x|\le\sqrt{\beta t}}\int_0^{t} G\left(t-s; y-x\right) \widetilde{U}^{(\beta, n)}_s(y) dB_s{(y)}. 
    \ees
    %Define the events
    %\[
    %    A_t(x) = \left( \wu_t(x) \le \epsilon \right ).
    %\]
    From Lemma \ref{log:wu:ratio} one can obtain that for some constant $\alpha(K)$ very large
    \[
        \log \bP \left( \wu_t(x) > K \right) \ge -\frac{1}{2}\log |x|
    \]
    for all $x\ge \alpha(K)$. For the same $x$
    \be \label{eq:u:le:ep}
        \bP \left( \wu_t(x) \le K \right) \le 1 -\frac{1}{\sqrt{|x|}}.
    \ee
    For the choice of $\beta=n$ and $k=2$ Lemma \ref{lem:U-Ubn} gives that
    \begin{align} \label{u_min_u^n}
    \begin{split}
        \bP\left( \wu_t(x) \le K, \wu^{(n, n)}_t(x) > 2K \right) &\le \bP\left( |\widetilde{U}_t(x) - \widetilde{U}_t^{(n, n)}(x)| > K \right) \\
        &\le \frac{C}{K^2} \left[\left(\frac{1}{2}\right)^n+ \exp\left(\frac{2\sqrt{n t}}{\mathcal R} \log\left(\frac{\mathcal Re\sqrt t}{\sqrt n}\right)\right)\right],
    \end{split}
    \end{align}
    for some constant $C$ depending on $t$ and $\|\wu_0\|_\infty$. Hence
    \begin{align}\label{eq:all:u:le:ep}
    \begin{split}
        &\bP \left( \bigcap_{x\in [n^4, 2n^4]\cap \Z} \wu_t(x) \le K \right) \\
        &\le \bP \left( \bigcap_{x\in [n^4, 2n^4]\cap \Z} \wu^{(n, n)}_t(x) \le 2K \right) + \frac{C}{K^2} \left[\frac{1}{2^n}+ \exp\left(\frac{2\sqrt{n t}}{\mathcal R} \log\left(\frac{\mathcal Re\sqrt t}{\sqrt n}\right)\right)\right]n^4.
    \end{split}
    \end{align}
    Define the integers $x_1=n^4$ and $x_{i+1} = n^4 + i\lceil 2n^{\frac{3}{2}}\sqrt{t}\rceil$ for every $i\ge1$. We will restrict attention to the integers $x_i \in [n^4, 2n^4]$, which are at most $\gamma_n = \frac{n^{4}}{2n^{3/2}\sqrt{t}}$ in number. Thanks to independence in Lemma \ref{lem:indep} and equation \eqref{u_min_u^n} we obtain
    \begin{flalign*}
        &\bP \left( \bigcap_{x\in [n^4, 2n^4]\cap \Z} \wu^{(n, n)}_t(x) \le 2K \right) &\\
        \quad&\le \bP \left( \bigcap_{i=1,\dots, \gamma_n} \wu^{(n, n)}_t(x_i) \le 2K \right) &\\
        \quad&=\prod_{i=1,\dots, \gamma_n} \left[ \bP \left( \wu^{(n, n)}_t(x_i) < 2K \right) \right] &\\
        \quad&\le \prod_{i=1,\dots, \gamma_n} \left[ \bP \left( \wu_t(x_i) \le 3K \right) + \bP \left( |\wu_t(x_i) - \wu^{(n, n)}_t(x_i)| > K \right) \right] &\\
        \quad&\le \left[ 1 - \frac{1}{\sqrt{2}n^2} + \frac{C}{K^2} \left\{\frac{1}{2^n}+ \exp\left(\frac{2\sqrt{n t}}{\mathcal R} \log\left(\frac{\mathcal Re\sqrt t}{\sqrt n}\right)\right)\right\} \right]^{\gamma_n},
    \end{flalign*}
    for a constant $C$, where in the second inequality we used equation \eqref{eq:u:le:ep} with $3K$ in place of $K$, as long as $n^4\ge \alpha(3K)$ is satisfied. Choose $n_1$ depending on $\RR, t$ and $K$ so that for all $n\ge n_1$
    \[
    \frac{Cn^4}{K^2} \left\{\frac{1}{2^{n}}+ \exp\left(\frac{2\sqrt{{n} t}}{\mathcal R} \log\left(\frac{\mathcal Re\sqrt t}{\sqrt {n}}\right)\right)\right\} < \frac{1}{2\sqrt{2}{n}^2}.
    \]
    If we choose $n\ge \alpha(3K)\vee n_1$ then equation \eqref{eq:all:u:le:ep} and the preceding calculation gives
    \begin{align*}
        \bP \left ( \sup_{x\in\Z} \wu_t(x) < K \right) &\le \bP \left( \bigcap_{x\in [n^4, 2n^4]\cap \Z} \wu_t(x) \le K \right) \\
        &\le \exp\left(-\frac{\sqrt{n}}{4\sqrt{2t}}\right) + \frac{1}{2\sqrt{2}{n}^2} \stackrel{n\to\infty}{\rightarrow} 0.
    \end{align*}
    As $K>0$ is arbitrary, this proves the Proposition.
\end{proof}

The following key proposition shows that the stochastic term of $\widetilde{U}$ is much smaller than the deterministic term when time is small.

\begin{proposition} \label{prop:det:noise} Fix any points $x, p\in \Z$ and $t>0$. For $\widetilde U$ with initial profile $\wu_0(\cdot) = M\cdot \mathbf{1}_{\{p\}}(\cdot)$,
\[ \bP\left(\left|\int_0^{t}\sum_{y\in \Z} G(t-s; y-x) \widetilde U_s(y) \, dB_s\left(y\right) \right| > \frac{M}{2} G(t; p-x)\right) \le 4  t\exp(t).\]
\end{proposition}
\begin{proof} 
We first bound $\E\left[\left(\widetilde U_t(y)\right)^2\right]$. Consider the Picard iterates defined as follows: $\widetilde U^{(\mathfrak{0})}_t(x)=M\cdot G(t; p-x)$, and for $\mathfrak n \ge 0$ let
\[ \widetilde U_t^{(\mathfrak{n+1})}(x) = M\cdot G(t; p-x) + \int_0^{t}\sum_{y\in  \Z} G(t-s; y-x) \widetilde U_s^{ (\mathfrak{n})}(y) \, dB_s\left(y\right).\]
With the definition 
\[\mathcal M_t^{(\mathfrak{n})}:= \sup_{y\in \Z} \E\left[\left(\frac{\widetilde U_t^{(\mathfrak{n})}(y)}{M \cdot G(t; p-y)}\right)^2\right], \]
it can be seen that 
\begin{align*} \mathcal M_t^{(\mathfrak{n+1})}\le 1+ \int_0^t \sum_{y \in  \Z}\left[\frac{G(t-s; p-y)\cdot G(s; y-x)}{G(t; p-x)}\right]^2 \mathcal{M}_s^{(\mathfrak{n})} ds.
\end{align*}
From the inequality $ \sum_i a_i^2 \le (\sum_i a_i)^2$ which holds for any sequence of non-negative real numbers $a_i$ and the Chapman-Kolmogorov equation we can deduce from the above display
\[ \mathcal M_t^{(\mathfrak{n+1})}\le 1+ \int_0^t\mathcal{M}_s^{(\mathfrak{n})} ds,\]
We know from equation (4.13) in \cite{geor-jose-khos-shiu} that the Picard iterates $ \wu_t^{(\mathfrak{n})}$ satisfy
\[
\sup_{\mathfrak{n}\ge 0}\; \sup_{y\in\Z} \E\left[\left|\widetilde U_t^{(\mathfrak{n})}(y)\right|^2\right] < \infty,
\]
and as a consequence $\sup_\mathfrak n \mathcal{M}_t^\mathfrak {(n)} < \infty$. Hence Gronwall's inequality gives
\[
    M_t^{(\mathfrak{n+1})} \le e^t.
\]
Therefore
\[ \E\left[\left(\widetilde U_t(x)\right)^2\right] \le M^2  \left(G(t; p-x)\right)^2 e^t.\]
We thus obtain
\begin{align*} \E\left[\left\{\int_0^{t}\sum_{y\in \Z} G(t-s; y-x) \widetilde U_s(y) \, dB_s\left(y\right)\right\}^2\right] &= \int_0^{t} \sum_{y \in\Z} \left(G(t-s; y-x)\right)^2 \E\left[\left(\widetilde U_s(y)\right)^2\right] ds \\
& \le M^2\int_0^{t}\sum_{y \in \Z} \left[G(t-s; y-x)\right]^2 \left[G(s; p-y)\right]^2e^s ds \\
&\le M^2\left[ G(t; p-x)\right]^2 t e^{t}
\end{align*}
The lemma follows from an application of Chebyshev's inequality. 
\end{proof}

\section{Proof of Theorem \ref{thm:main}}\label{sec:thm}
%\normalfont
Fix $\delta>0$ and $L>0$ and a point $x\in\Z$. We will show $\bP(\sup_{t\le 3\delta} U_t(x) \ge L) =1$. Increasing $L$ to infinity will show $\bP(\sup_{t\le 3\delta} U_t(x) =\infty)=1$, and then the arbitrariness of $\delta$ and $x$ will establish that $U$ explodes instantaneously everywhere.

Choose and fix $n_0$ and $K_b$ as in \eqref{eq:kb}. Recall from \eqref{eq:z:sde} the SDE
\[ dZ_t= \frac{1}{n_0} b(Z_t) \, dt + Z_t\, dB_t(x).\]
%Let $Z_t^{(J)}$ solve the same SDE as $Z_t$ but with truncated drift $b(\cdot \wedge J)$ instead. A drift comparison gives that $Z^{(J)}\le Z^{(J+1)} \le \dots \le Z$. 
Using Proposition \ref{prop:explosion} find $K>2K_b$ such that 
\be \label{eq:Z:explode} \bP_K\left(Z_t \text{ explodes by time } \delta\right) \ge 1-\delta.\ee

Next consider the process $U$ between time $0$ and $\delta$. Clearly this is {\it bounded below} by the process $\widetilde U$ without drift given in \eqref{eq:ISDE:2}.  Let $p$ be the closest {\it random} point to $x$ such that $\widetilde U_{\delta}(p) \ge K$. Such a point exists with probability $1$ due to Proposition \ref{cor:sp:growth}. By Mueller's comparison principle the profile of $U$ at time $\delta$ is bounded below by $K \mathbf{1}_{\{p\}}(\cdot)$ for the random point $p$.

Together with the fixed random point $p$ let us find a random level $M=M(p)$ depending on $p$ such that
\be\label{eq:M:choice} M \ge \frac{2L}{G(\delta; p-x)}. \ee
Consider the process $\overline U^{(2M)}_t(z), \, t\ge \delta, \, z\in \Z$ solving 

\[ 
d \overline U^{(2M)}_t(z) = \left(\mathscr{L} \overline U^{(2M)}_t\right)(z)\, dt + b\left(\overline U^{(2M)}_t(z) \wedge 2M\right)\, dt + \overline U^{(2M)}_t(z)\, dB_t(z),
\]
with {\it random} initial profile $\overline U^{(2M)}_{\delta}(z) = K \mathbf{1}_{\{p\}}(\cdot)$. Clearly
\be \label{eq:U:oU} \bP\left(\sup_{t\in[0, 3\delta]} U_t(x) \ge L\right) \ge \bP\left(\sup_{t\in [\delta, 3\delta] } \overline U^{(2M)}_t(x) \ge L \right).\ee

Now let $V^{(n,2M)}$ be the Alternating process corresponding to the ISDE $\overline{U}^{(2M)}$. Let $Z_t^{(2M)}$ solve the same SDE as $Z_t$ but with $b$ replaced by $b(\cdot \wedge 2M)$ instead and driven by the Brownian motion at random site $p$. An appeal to Corollary \ref{cor:UV} gives
\begin{align}\label{eq:step1}
\begin{split}
    \bP\left(\sup_{t\in [\delta, 2\delta]} \overline U_t^{(2M)}(p) \ge M \right) &\ge \limsup_{n\to \infty} \bP\left(\sup_{t\in [\delta, 2\delta]} V_t^{(n,2M)}(p)\ge M \right) \\
    &\ge \bP\left(\sup_{t\in [\delta, 2\delta]} Z_t^{(2M)}\ge M \right) \\
    &\ge 1- \delta,
\end{split}
\end{align}
where second line is due to our choice of $K>2K_b$ and the following consequence of Proposition \ref{prop:compare}: before the level $2M$ is attained for the first time by $Z_t^{2M}$, we have $Z_t=Z_t^{(2M)}$ and $V_t^{(n,2M)}(x)$ has to be at least as large as $Z_t$. The third line is because of \eqref{eq:Z:explode}.

Let $\tau$  be the random time defined as
\[ \tau := \inf\left\{t\ge \delta: \,\overline U^{(2M)}_t(p) \ge M \right\},\]
and equal to $\infty$ if the set is empty. The equation \eqref{eq:step1} is equivalent to $\bP(\tau \le 2\delta) \ge 1-\delta$. On the event $\tau <\infty$, consider the process $\widehat U_t(z), \, t\ge \tau, \, z\in \Z$ solving 
\[ d \widehat U_t(z) = \left(\mathscr{L} \widehat U_t\right)(z)\, dt +   \widehat U_t(z)\, dB_t(z),\]
with {\it random} initial profile $\widehat U_{\tau}(z) = M\mathbf{1}_{\{p\}}(\cdot)$. Therefore %We have by the strong Markov property and Mueller's comparison principle 
\be \label{eq:step2}
\begin{split}
&\bP\left(\sup_{t\in [\delta, 3\delta] }  \overline U^{(2M)}_t(x) \ge L\right)\\
&\ge \bP\left(\sup_{t\in [\delta, 3\delta] }  \overline U^{(2M)}_t(x) \ge L,\, \sup_{t\in [\delta, 2\delta] }  \overline U^{(2M)}_t(p) \ge M\right)  \\
& = \bP\left(\sup_{t\in [\delta, 2\delta] }  \overline U^{(2M)}_t(p) \ge M\right) \cdot \bP\left(\sup_{t\in [\delta, 3\delta] }  \overline U^{(2M)}_t(x) \ge L \Big \vert\sup_{t\in [\delta, 2\delta] }  \overline U^{(2M)}_t(p) \ge M\right)  \\
&\ge (1-\delta)\cdot \bP\left(\sup_{t\in [\tau, \tau+\delta]} \widehat U_t(x) \ge L\right),
\end{split}
\ee
the last line following from the strong Markov property and Mueller's comparison principle. Moreover, one sees from Proposition \ref{prop:det:noise} and the choice of $M$ in \eqref{eq:M:choice} that
\be \label{eq:step3} \bP\left(\sup_{t\in [\tau, \tau+\delta]} \widehat U_t(x) \ge L\right) \ge 1- 4\delta e^{\delta}.\ee
Putting together \eqref{eq:U:oU}, 
\eqref{eq:step2} and \eqref{eq:step3} we obtain
\bean
%\begin{split}
\bP\left(\sup_{t\in [0,3\delta]} U_t(x) \ge L\right) &\ge (1-\delta)(1-4\delta e^{\delta}).
%\end{split}
\eean
Finally observe that the events $\{\sup_{t\le 3\delta} U_t(x) \ge L\}$ are decreasing as $\delta$ goes to $0$, while the lower bound increases to $1$ as $\delta$ goes to $0$. Therefore 
\[  \bP\left(\sup_{t\in [0, 3\delta]} U_t(x) \ge L\right) =1,\]
completing the proof of the theorem. \qed

\subsection*{Acknowledgements}
We would like to thank the two anonymous referees for their careful reading of the manuscript and their many suggestions which have greatly improved this paper.

%\bibliography{blowup}
%\bibliographystyle{abbrv}

\end{document}